\documentclass[twocolumn,amsthm]{autart}

\usepackage{graphicx}

\usepackage{amsmath}
\usepackage{url}
\usepackage{tikz}
\usepackage{pgfplots}
\usepackage{soul}
\usepackage{tabularx}
\usepackage{arydshln}
\usepackage{siunitx}

\pgfdeclarehorizontalshading{parula}{140bp}{
    rgb(0bp)=(0.242200,0.150400,0.660300);
    rgb(20bp)=(0.242200,0.150400,0.660300);
    rgb(30bp)=(0.279100,0.262550,0.906450);
    rgb(40bp)=(0.264700,0.403000,0.993500);
    rgb(50bp)=(0.176600,0.547800,0.953800);
    rgb(60bp)=(0.108500,0.666900,0.873400);
    rgb(70bp)=(0.070400,0.745700,0.725800);
    rgb(80bp)=(0.280900,0.796400,0.526600);
    rgb(90bp)=(0.622900,0.786850,0.257150);
    rgb(100bp)=(0.918400,0.730800,0.189000);
    rgb(110bp)=(0.982100,0.831100,0.180200);
    rgb(120bp)=(0.976900,0.983900,0.080500);
    rgb(140bp)=(0.976900,0.983900,0.080500)
}
\pgfplotsset{compat=1.18}

\usepackage{amssymb}
\usepackage{amsfonts}
\usepackage{mathtools}
\usetikzlibrary{shadings}
\usepackage{bm}
\usepackage[inline]{enumitem}
\usepackage{xcolor}
\usepackage{subcaption}

\DeclareMathOperator{\rank}{rank}
\DeclareMathOperator{\tr}{tr}
\newcommand{\hop}{\mathsf{H} }

\newcounter{assumption}
\newcounter{theorem}

\newcounter{definition}

\newcounter{lemma}

\allowdisplaybreaks

\makeatletter
\DeclareFontFamily{OMX}{MnSymbolE}{}
\DeclareSymbolFont{MnLargeSymbols}{OMX}{MnSymbolE}{m}{n}
\SetSymbolFont{MnLargeSymbols}{bold}{OMX}{MnSymbolE}{b}{n}
\DeclareFontShape{OMX}{MnSymbolE}{m}{n}{
    <-6>  MnSymbolE5
   <6-7>  MnSymbolE6
   <7-8>  MnSymbolE7
   <8-9>  MnSymbolE8
   <9-10> MnSymbolE9
  <10-12> MnSymbolE10
  <12->   MnSymbolE12
}{}
\DeclareFontShape{OMX}{MnSymbolE}{b}{n}{
    <-6>  MnSymbolE-Bold5
   <6-7>  MnSymbolE-Bold6
   <7-8>  MnSymbolE-Bold7
   <8-9>  MnSymbolE-Bold8
   <9-10> MnSymbolE-Bold9
  <10-12> MnSymbolE-Bold10
  <12->   MnSymbolE-Bold12
}{}

\let\llangle\@undefined
\let\rrangle\@undefined
\DeclareMathDelimiter{\llangle}{\mathopen}%
                     {MnLargeSymbols}{'164}{MnLargeSymbols}{'164}
\DeclareMathDelimiter{\rrangle}{\mathclose}%
                     {MnLargeSymbols}{'171}{MnLargeSymbols}{'171}
\makeatother

\begin{document}

\begin{frontmatter}

\title{Moving-Boundary~Port-Hamiltonian~Systems: Mathematical~Formalization, Analysis, and Simulation}

\thanks[footnoteinfo]{Corresponding author: T.~J.~Meijer.}

\author[TUe_ME]{Tomas~J.~{Meijer}}\ead{t.j.meijer@tue.nl},
\author[TUe_EE]{Amritam~{Das}}\ead{am.das@tue.nl},
\author[TUe_EE]{Siep~{Weiland}}\ead{s.weiland@tue.nl}

\address[TUe_ME]{Department of Mechanical Engineering, Eindhoven University of Technology, 5600 MB Eindhoven, The Netherlands}
\address[TUe_EE]{Department of Electrical Engineering, Eindhoven University of Technology, 5600 MB Eindhoven, The Netherlands}

\begin{abstract}
    In this paper, we consider linear boundary port-Hamiltonian distributed parameter systems on a time-varying spatial domain. We derive the specific time-varying Dirac structure that these systems give rise to and use it to formally establish a new class of moving-boundary port-Hamiltonian systems by showing that these distributed parameter systems on a time-varying spatial domain admit a port-Hamiltonian representation. We demonstrate that our results can be leveraged to develop a spatial discretization scheme with dynamic meshing for approximating the telegrapher's equations on a time-varying spatial domain, which we subsequently verify numerically.
\end{abstract}

        
\end{frontmatter}

\section{Introduction}
\label{sec:introduction}
Boundary port-Hamiltonian (pH) distributed parameter systems (DPSs) provide an extension of classical Hamiltonian dynamics to open systems that admit energy flow through the boundary of their spatial domain~\cite{vanderSchaft2002,Duindam2014}. This extension is realized by introducing so-called port variables on the boundary of the spatial configuration space through which power flow across the boundary into (or out of) the system is realized. Such a boundary pH system admits a Dirac structure, or more precisely a Stokes-Dirac structure, derived from the Hamiltonian differential operator~\cite{Duindam2014}. Importantly, their specific mathematical structure and the resulting adherence to physical conservation laws (e.g., conservation of energy, momentum, et cetera) have proven instrumental for passivity- and energy-based (boundary) control of DPSs~\cite{Hamroun2010,Macchelli2018,Augner2020,Jacob2012,Califano2022}. PH systems have also given rise to so-called structure-preserving discretization, see, e.g.,~\cite{Golo2004,Mehrmann2019,Brugnoli2022,Voss2011,Wang2017}, and model-order reduction methods, see, e.g.,~\cite{Gugercin2012,Schwerdtner2023,Breiten2022}, which aim to preserve their mathematical structure, desirable properties and physical interpretability so that they can be exploited for analysis or control purposes~\cite{Lou2021}.

Although pH-DPSs are commonly defined on a fixed spatial domain, many classical problems and modern applications require considering a time-varying spatial domain instead. An example of such a classical problem is the Stefan problem, which describes phase transitions in matter and features a moving interface between the different phases~\cite{Rubinshtein1971}. Time-varying spatial domains are also useful in mechanical systems that structurally deform over time, see, e.g.,~\cite{Rees2012}. A modern application, in which such thermomechanical deformations are highly relevant~\cite{vandenHurk2020}, is found in lithography equipment for the production of integrated circuits (ICs), see, e.g.,~\cite{Butler2011,vandenHurk2020,Meijer2024-phd-thesis,Veldman2020}. In each of these problems/applications, the governing differential equations are derived from physical conservation laws. However, formal characterization of these laws, thereby, preservation of them, is generally ignored when using (approximations of) these models for controller design, see, e.g.,~\cite{vandenHurk2020}, or simulation, see, e.g.,~\cite{Gibou2005}. In a passivity-based or energy-based control context, this is detrimental because stability properties or performance guarantees do not apply to the underlying system.

\begin{figure}[!t]
    \centering    
        \begin{subfigure}{.95\linewidth}
            \includegraphics[width=\linewidth]{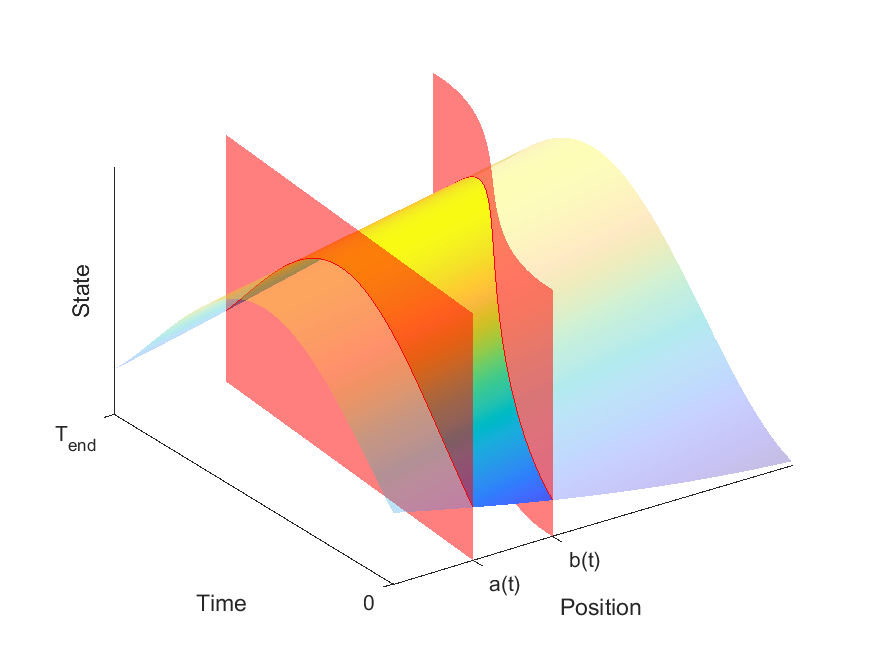}
            \caption{Side view.}
        \end{subfigure} 
        \begin{subfigure}{.95\linewidth}
            \includegraphics[width=\linewidth]{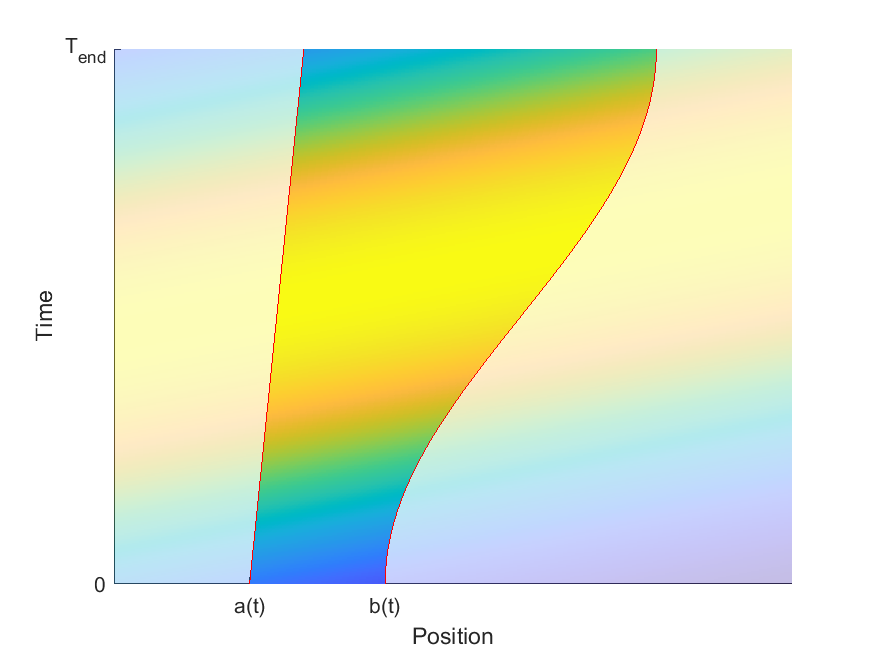}
            \caption{Top view.}
        \end{subfigure}  
    \caption{Solution to the transport equation inside (\protect\tikz\protect\shade[shading=parula,shading angle=135,opacity=1] (0,0) rectangle (0.2,0.2);) and outside (\protect\tikz\protect\shade[shading=parula,shading angle=135,opacity=0.3] (0,0) rectangle (0.2,0.2);) a spatial domain with moving boundaries (\protect\tikz\protect\filldraw[fill=red,opacity=0.5,draw=none] (0,0) rectangle (0.2,0.2);).}
    \label{fig:transport-eq-tv-bnds}
\end{figure}
A port-Hamiltonian formulation of DPSs on time-varying spatial domains can potentially circumvent the aforementioned issues by providing a framework in which the different power flows and conservation laws can be systematically characterized and analyzed. Figure~\ref{fig:transport-eq-tv-bnds} depicts such a time-varying spatial domain, which we will denote by $\mathbb{S}_{ab}(t)\coloneqq [a(t),b(t)]$, for the transport equation along with its solution. Importantly, the physical quantities typically extend beyond the system's time-varying spatial domain, as can be seen in Fig.~\ref{fig:transport-eq-tv-bnds}. From the perspective of $\mathbb{S}_{ab}(t)$, the movement of the boundary results in an additional flow of physical quantities into or out of $\mathbb{S}_{ab}(t)$. Hence, additional dynamics are introduced that depend on the velocity of the boundaries. As we will see, the movement of the boundaries also introduces additional power flows into or out of the port-Hamiltonian system. Existing works consider port-Hamiltonian systems separated by a moving interface, see, e.g.,~\cite{Diagne2013,Vincent2020,Kilian2023}, which are a special case of such moving-boundary pH systems. We formalize these insights in a novel extension to boundary pH-DPSs, which we call \emph{moving-boundary port-Hamiltonian systems}. Facilitation of the simulation and control of these systems requires to also develop spatial discretization schemes for them. Naturally, when considering time-varying spatial domains, the mesh used to discretize this domain also needs to be dynamic such that it can be adapted, and, at all times, appropriately reflects the spatial domain of the system. Although there exist methods that iteratively refine the mesh based on the obtained residuals, see, e.g.,~\cite{Carey1981}, to the best of the authors' knowledge, real-time dynamic meshing remains an open problem. A key challenge in developing such a scheme is to incorporate the dynamic boundaries between neighboring elements while adhering to the underlying physical conservation laws of the system. Moving-boundary pH systems  promise to provide an important step toward an intuitive and modular solution to this challenge by considering each element in a dynamic grid as an individual instance of such a system. Dynamic meshing is also powerful for applications involving DPSs with static spatial domains but with, for instance, moving actuators or sensors, see, e.g.,~\cite{Veldman2020}, or control tasks that focus on a time-varying region or point of interest in the spatial domain, see, e.g.,~\cite{vandenHurk2020}. In such applications, the use of a dynamic grid to subdivide the system into finite elements enables us to maintain a high spatial resolution in the vicinity of these moving sensors, actuators, and regions or points of interest at all times. 

This paper aims to answer the following open questions: Firstly, how do we formalize the underlying pH (i.e., Dirac) structure for pH-DPSs on time-varying spatial domains? Secondly, how do we leverage this structure in discretization schemes for such systems so that the discretization can adapt to their varying spatial domain? 

To answer these questions, we present the following important contributions. Firstly, we model the dynamics of pH-DPSs on one-dimensional time-varying spatial domains, study their conservation of energy, and derive the corresponding power balance equations. We show that these systems give rise to a specific \emph{time-varying} Stokes-Dirac structure, which we use to formally establish the novel class of moving-boundary port-Hamiltonian systems. Finally, we demonstrate, for a time-varying segment of a lossless transmission line (TL), that our formulation of moving-boundary pH systems in terms of a specific Stokes-Dirac structure naturally yields a spatial discretization scheme with a dynamic grid, which we can use for simulation or control purposes and for which we showcase numerical results. 

The remainder of the paper is organized as follows. Section~\ref{sec:prelim} introduces relevant notation and preliminaries regarding pH-DPSs. In Section~\ref{sec:main-results}, we present our main results, which consist of a thorough analysis and formalization of moving-boundary pH systems and the underlying Dirac structure. We propose a spatial discretization of the TL on $\mathbb{S}_{ab}(t)$ and use it to simulate results in Section~\ref{sec:numerical-results}. Finally, Section~\ref{sec:conclusions} provides some conclusions, and all proofs can be found in the Appendix.

\section{Notation and preliminaries}\label{sec:prelim}
The set of real numbers, non-negative real numbers, positive real numbers, and non-negative natural numbers are denoted, respectively, by $\mathbb{R}=(-\infty,\infty)$, $\mathbb{R}_{\geqslant 0}=[0,\infty)$, $\mathbb{R}_{>0}=(0,\infty)$, and $\mathbb{N}=\{0,1,2,\hdots\}$. We use $I_n$ to denote the $n$-by-$n$ identity matrix. For a vector $v\in\mathbb{C}^{n}$, we use $v^\top$ and $v^\hop$ to denote, respectively, the transpose and complex-conjugate transpose of $v$. The sets of square-integrable functions and $k$-times continuously-differentiable functions on $X$ with co-domain $Y$ are denoted by, respectively, $L_2(X,Y)$ and $C^k(X,Y)$. For the latter, the order of differentiation is omitted in the case $k=1$. We use $\partial_x$ and $\delta_x$, respectively, to denote the partial and variational derivative with respect to $x$, and, for $f\in C(\mathbb{R},\mathbb{R}^{n})$, we denote $\dot{f}=(\mathrm{d}f)/(\mathrm{d}t)$. Let $\mathbb{K}$ denote a field that is either the real numbers $\mathbb{R}$ or the complex numbers $\mathbb{C}$, and let $\left\langle\cdot,\cdot\right\rangle_{L_2}$ and $\left\langle\cdot,\cdot\right\rangle_{\mathbb{K}^n}$ denote, respectively, the Hilbert space inner product, i.e., $\left\langle u,v\right\rangle =\int_a^b u^\top(s)v(s)\,\mathrm{d}s$ for any $u,v\in L_2([a,b],\mathbb{R}^{n})$, and the Hermitian inner product on $\mathbb{K}^{n}$, i.e., $\left\langle u,v\right\rangle_{\mathbb{K}^n} = u^\hop v$ for any $u,v\in\mathbb{K}^n$. Let $(u_1,u_2,\hdots,u_q)\coloneqq \operatorname{col}(u_1,u_2,\hdots,u_q)$ for any vectors $u_i\in\mathbb{K}^{n_i}$, $i\in\{1,2,\hdots,q\}$. We denote the trace operator $\tr~:~L_2([a,b],\mathbb{R}^n)\rightarrow\mathbb{R}^{2n}$, i.e., $\tr(x)=(x(a),x(b))$, $x\in L_2([a,b],\mathbb{R}^{n})$. 

\subsection{Linear first-order boundary pH-DPSs}
Before we introduce the considered class of systems, we first provide a definition of a Dirac structure below. 
\setcounter{thm}{\thedefinition}\stepcounter{definition}
\begin{defn}\label{def:Dirac}
    Let $\mathcal{F}$ be a vector space of flows over the field $\mathbb{K}$ endowed with a non-degenerate bilineair form $\left\langle \cdot\middle|\cdot\right\rangle~\colon~\mathcal{F}\times\mathcal{E}\rightarrow\mathbb{K}$. On the product space $\mathcal{B}=\mathcal{F}\times\mathcal{E}$, where the space of efforts $\mathcal{E}$ is the dual space of $\mathcal{F}$, i.e., $\mathcal{E}=\mathcal{F}^*$, consider the symmetric bilinear form $\left\llangle \cdot,\cdot\right\rrangle~\colon~\mathcal{B}\times\mathcal{B}\rightarrow\mathbb{K}$ given by $\left\llangle (f_1,e_1),(f_2,e_2)\right\rrangle\coloneqq \left\langle f_2\middle|e_1\right\rangle+\left\langle f_1\middle|e_2\right\rangle$. A \emph{Dirac structure} is a subspace $\mathcal{D}\subset\mathcal{B}$ such that $\mathcal{D}=\mathcal{D}^\perp$, where $\mathcal{D}^\perp\coloneqq \{d_1\in\mathcal{B}\mid \left\llangle d,d_1\right\rrangle=0 \text{ for all }d\in\mathcal{D}\}$ is the orthogonal complement with respect to $\left\llangle\cdot,\cdot\right\rrangle$.
\end{defn}
Typically, Dirac structures are defined over the field of reals, see, e.g.,~\cite{Courant1988,vanderSchaft2002}, however, the definition naturally carries over to linear spaces over an arbitrary field~\cite{Jeltsema2015}, which yields the above definition. As we will see, the resulting additional generality is crucial in formulating the underlying Dirac structure for the class of moving-boundary port-Hamiltonian systems studied in this paper. On a related note, complex-valued port-Hamiltonian systems are also considered in, e.g.,\cite{Beattie2018}.

Next, we introduce linear first-order boundary pH-DPSs. Consider the first-order Hamiltonian operator 
\begin{equation*}
    \mathcal{J}\coloneqq J_0+J_1\partial_s,
\end{equation*}
with $J_0=-J_0^\top\in\mathbb{R}^{n\times n}$ and $J_1=J_1^\top\in\mathbb{R}^{n\times n}$. This Hamiltonian operator is associated with a specific type of Dirac structure, called \emph{Stokes-Dirac structure}.
\setcounter{thm}{\thelemma}\stepcounter{lemma}
\begin{lem}{\cite[Thm.~2.7]{Villegas2007}}\label{lem:stokes-dirac}
    Let $\rank J_1=r$, $\mathcal{F}=L_2([a,b],\mathbb{R}^{n})\times \mathbb{R}^{r}$, $\mathcal{E}=\mathcal{F}^*$, and let $\mathcal{B}=\mathcal{F}\times\mathcal{E}$ be equipped with the non-degenerate bilinear form $\left\langle\cdot\middle|\cdot\right\rangle_\mathcal{B}~\colon~\mathcal{B}\rightarrow\mathbb{R}$ given by
    \begin{equation*}
        \left\langle(f,f_{\partial})\middle|(e,e_{\partial})\right\rangle_\mathcal{B}= \left\langle f,e\right\rangle_{L_2} - \left\langle f_\partial,e_{\partial}\right\rangle_{\mathbb{R}^r}.
    \end{equation*}
    Then, there exist a symmetric matrix $S_1\in\mathbb{R}^{n\times n}$ and a full rank matrix $M\in\mathbb{R}^{r\times n}$ such that the vector subspace
    \begin{align*}\label{eq:DJ}
        \mathcal{D}_{\mathcal{J}}&\coloneqq\left\{\left(f,f_{\partial},e,e_{\partial}\right)\in\mathcal{B}\,\middle|\,e\in\mathcal{E},\mathcal{J}e\in\mathcal{F},\right.\\
        &\qquad \left.f=\mathcal{J}e,~\begin{bmatrix}
            f_{\partial}\\
            e_{\partial}
        \end{bmatrix} = \begin{bmatrix}
            S_1 & -S_1\\
            I_n & I_n
        \end{bmatrix}\tr(Me)\right\},
    \end{align*}
    is a Dirac structure with respect to the symmetric bilinear form $\left\llangle \cdot,\cdot\right\rrangle_{\mathcal{B}}~\colon~\mathcal{B}\times\mathcal{B}\rightarrow\mathbb{R}$ given by
    \begin{align*}
        &\left\llangle \left(f^1,f^1_{\partial},e^1,e^1_{\partial}\right),\left(f^2,f^2_{\partial},e^2,e^2_{\partial}\right)\right\rrangle_{\mathcal{B}} = \\
        &\quad\left\langle \left(f^1,f^1_{\partial}\right)\middle|\left(e^2,e^2_{\partial}\right)\right\rangle_{\mathcal{B}} + \left\langle \left(f^2,f^2_{\partial}\right)\middle|\left(e^1,e^1_{\partial}\right)\right\rangle_{\mathcal{B}}.
    \end{align*}		
\end{lem}
In addition to the effort $e$ and flow $f$, the Stokes-Dirac structure $\mathcal{D}_{\mathcal{J}}$ introduced in Lemma~\ref{lem:stokes-dirac} incorporates boundary ports $(f_\partial,e_\partial)$ through which the system interacts with its environment, i.e., through which the boundary conditions of the system are imposed. To construct these boundary ports we need to construct the matrices $S_1$ and $M$. We refer to~\cite{Villegas2007} for details on how to construct $S_1$ and $M$ in the general case, however, if $J_1$ is full rank, $M=I_n$ and $S_1=-\frac{1}{2}J_1$~\cite{LeGorrec2005}. 

Using $\mathcal{D}_{\mathcal{J}}$, we now formally introduce the linear first-order boundary pH-DPSs considered in this paper.
\setcounter{thm}{\thedefinition}\stepcounter{definition}
\begin{defn}{\cite{Maschke2023,Villegas2007}}\label{def:lin-fo-bound-pH-DPSs}
    The linear boundary pH system on the state space $\mathcal{X}= L_2([a,b],\mathbb{R}^n)$ with respect to the Stokes-Dirac structure $\mathcal{D}_{\mathcal{J}}$ associated with $\mathcal{J}$ and generated by the \emph{Hamiltonian} functional $H~\colon~\mathcal{X}\rightarrow\mathbb{R}_{\geqslant 0}$ given by 
    \begin{equation*}
        H(x) = \frac{1}{2}\int_a^b x^\top(s) Qx(s)\,\mathrm{d}s,
    \end{equation*} 
    where $Q=Q^\top\in\mathbb{R}^{n\times n}$ is a positive definite matrix, is the dynamical system defined by
    \begin{equation}
        \left(\partial_t x,f_{\partial},Qx,e_{\partial}\right)\in\mathcal{D}_{\mathcal{J}},\label{eq:pH-DPS}
    \end{equation}
    where $x\in \mathcal{X}$ is its state.
\end{defn}
Due to its physical interpretation, the state $x$ in Definition~\ref{def:lin-fo-bound-pH-DPSs} is referred to as the energy variable of the system, whereas $e\coloneqq \delta_x H = Qx$ in~\eqref{eq:pH-DPS} is called the co-energy variable. These systems arise naturally in first-principle modelling of physical systems, see, e.g.,~\cite{vanderSchaft2002,Villegas2007,Voss2011,Hamroun2006,Baaiu2009,LeGorrec2005}, and have proven to be useful in passivity-based boundary control of PDEs~\cite{Jacob2012}. 

Note that, since $Q$ is a constant matrix, the scope of Definition~\ref{def:lin-fo-bound-pH-DPSs} only covers homogeneous PDEs. While the results presented in the results presented in the sequel can be naturally extended to accommodate non-homogeneous PDEs, doing so introduces many additional terms and, thus, has been avoided in the present paper for the sake of compactness.

\section{Moving-boundary port-Hamiltonian systems}\label{sec:main-results}
\subsection{Dynamics of moving-boundary pH systems}
We consider linear first-order boundary pH-DPSs~\eqref{eq:pH-DPS} on a time-varying spatial domain $\mathbb{S}_{ab}(t)=\left[a(t),b(t)\right]$ with $a(t)<b(t)$, $t\in\mathbb{T}$ (refer to Fig.~\ref{fig:transport-eq-tv-bnds} for an illustration). This leads us to define the time-dependent state space $\mathcal{X}~\colon~\mathbb{T}\rightarrow L_2([a(\cdot),b(\cdot)],\mathbb{R}^{n})$, defined at time $t\in\mathbb{T}$ as $\mathcal{X}(t)=L_2([a(t),b(t)],\mathbb{R}^{n})$, and the Hamiltonian $H~\colon~\mathcal{X}\times\mathbb{T}\rightarrow\mathbb{R}_{\geqslant 0}$ as
\begin{equation*}
        H(x,t) = \frac{1}{2}\int_{a(t)}^{b(t)} x^\top(s) Qx(s)\,\mathrm{d}s,
\end{equation*}
with the state $x\in\mathcal{X}(t)$ for all $t\in\mathbb{T}$, which incorporates the time-varying nature of the spatial domain $\mathbb{S}_{ab}(t)=\left[a(t),b(t)\right]$ in the integration bounds. Similarly, in Lemma~\ref{lem:stokes-dirac}, this causes us to introduce the time-varying space of flows $\mathcal{F}~\colon~\mathbb{T}\rightarrow L_2([a(\cdot),b(\cdot)],\mathbb{R}^{n})\times \mathbb{R}^{r}$, whereby the space of efforts $\mathcal{E}(t)=\mathcal{F}^*(t)$, the product space $\mathcal{B}(t)=\mathcal{F}(t)\times\mathcal{E}(t)$, and the Dirac structure $\mathcal{D}_{\mathcal{J}}(t)\subset\mathcal{B}(t)$, $t\in\mathbb{T}$, also become time-dependent. 

Regarding $\mathbb{S}_{ab}(t)$, we assume the following.
\setcounter{thm}{\theassumption}\stepcounter{assumption}
\begin{assum}\label{asm:bounds}
    The boundaries of the spatial domain $\mathbb{S}_{ab}(t)$ satisfy the following:
    \begin{enumerate}[label=(\roman*)]
        \item $a,b\in C(\mathbb{T},\mathbb{R})$;
        \item $a(t)<b(t)$ for all $t\in\mathbb{T}$;
        \item $\dot{a}(t)\dot{b}(t)\geqslant 0$ for all $t\in\mathbb{T}$.\label{item:same-sign}
    \end{enumerate}
\end{assum}
Assumption~\ref{asm:bounds} means that $\mathbb{S}_{ab}(t)$ evolves over time in a continuous manner with $a(t)$ and $b(t)$ retaining their order for all $t\in\mathbb{T}$. Moreover, item \ref{item:same-sign}, which is needed for technical reasons in proving our results, ensures that both boundaries always move in the same direction.

 To deal with the time-varying integration bounds, we introduce the following time-varying change of spatial variable $s(t, \hat{s})$, for all $\hat{s}\in[0,1]$ and $t\in\mathbb{T}$,
\begin{align}
    s(t, \hat{s}) =& h(t,\hat{s}) \notag\\:=& a(t) + \left(b(t)-a(t)\right)\hat{s},\label{eq:coc}  
\end{align}
This change of variable allows the Hamiltonian $H$ to be expressed as
\begin{equation}\label{eq:H-coc}
    H(x,t) = \frac{1}{2}\int_0^1 (b(t)-a(t))x^\top(h(t,\hat{s}))Qx(h(t,\hat{s}))\,\mathrm{d}\hat{s},
\end{equation}
 where we used the fact that $\mathrm{d}s(t,\hat{s}) = (b(t)-a(t))\mathrm{d}\hat{s}$. Moreover,~\eqref{eq:coc} also leads to a redefined state space $\hat{\mathcal{X}}~:~\mathbb{T}\rightarrow L_2([0,1],\mathbb{R}^n)$, and the new state $\hat{x}\in\hat{\mathcal{X}}(t)$ for all $t\in\mathbb{T}$, which is related to $x$ by
\begin{equation}\label{eq:xhat}
    \hat{x}(t,\hat{s})\coloneqq \sqrt{b(t)-a(t)}x(t,h(t,\hat{s})).
\end{equation}
Here and in the sequel, we abuse notation to write $x(t,s)$ and $\hat{x}(t,\hat{s})$ if, respectively, $x(t,\cdot)\in\mathcal{X}(t)$ and $\hat{x}(t,\cdot)\in\hat{\mathcal{X}}(t)$ for all $t\in\mathbb{T}$. 

Using the time-varying state transformation in~\eqref{eq:xhat}, we can re-express~\eqref{eq:H-coc} in a form that is similar to $H$ in Definition~\ref{def:lin-fo-bound-pH-DPSs}. To this end, we introduce a new Hamiltonian $\hat{H}~\colon~\hat{\mathcal{X}}\times \mathbb{T}\rightarrow\mathbb{R}_{\geqslant 0}$, given by
\begin{equation}\label{eq:Hhat}
    \hat{H}(\hat{x},t) = \frac{1}{2}\int_0^1\hat{x}^\top(\hat{s})Q\hat{x}(\hat{s})\,\mathrm{d}\hat{s},
\end{equation}
with $\hat{x}\in\hat{\mathcal{X}}(t)$. Importantly, due to~\eqref{eq:xhat}, $\hat{H}$ satisfies
\begin{equation}
    H(x,t) = \hat{H}(\hat{x},t)\text{ for all }t\in\mathbb{T}.
\end{equation}

In Lemma~\ref{lem:dyn}, we present the dynamics that govern the new energy variables. To this end, let $\hat{e}\in\hat{\mathcal{E}}(t)=\hat{\mathcal{X}}^*(t)$, with $\hat{\mathcal{E}}~\colon~\mathbb{T}\rightarrow L_2([0,1],\mathbb{R}^{n})$, be the new co-energy variables defined as 
\begin{equation}\label{eq:efforts}
    \hat{e}(t,\hat{s}) \coloneqq \delta_{\hat{x}}\hat{H}(\hat{x},t) = Q\hat{x}(t,\hat{s}).
\end{equation}
\begin{lem}\label{lem:dyn}
    The dynamics of $\hat{x}$ are governed by 
    \begin{align}
        \partial_t\hat{x}(t,\hat{s}) &= \hat{\mathcal{J}}(t)\hat{e}(t,\hat{s})-\frac{1}{2}\frac{\dot{b}(t)-\dot{a}(t)}{b(t)-a(t)}\hat{x}(t,\hat{s}) +\label{eq:dyn}\\
        &\quad  \partial_{\hat{s}}\left(\frac{\dot{a}+(\dot{b}-\dot{a})\hat{s}}{b-a}\hat{x}(t,\hat{s})\right),\nonumber
    \end{align}
    where $\hat{\mathcal{J}}(t)\coloneqq J_0+(1/(b(t)-a(t)))J_1\partial_{\hat{s}}$.    
\end{lem}
\begin{proof}
    See Appendix~\ref{app:dynamics} for a detailed proof.
\end{proof}
By inspection of~\eqref{eq:dyn}, we note that, aside from the dynamics of the original system~\eqref{eq:pH-DPS} captured by the $\hat{\mathcal{J}}(t)\hat{e}(t,\hat{s})$-term, the new state dynamics contain two additional terms due to the movement of the boundaries. Interestingly, each of these terms admits a physically insightful interpretation: 
\begin{enumerate}[label=(\roman*)]
    \item The term $-(\dot{b}(t)-\dot{a}(t))/(2(b(t)-a(t))))\hat{x}$ corresponds to the net {\bf(de)compression} of the spatial domain $\mathbb{S}_{ab}(t)$, and the state $x(t,s)$ contained therein, as a result of the moving boundaries;
    \item Any point of interest $\hat{s}\in[0,1]$ maps to a point $s(t,\hat{s})\in[a(t),b(t)]$ through~\eqref{eq:coc}. Hence, the $\partial_{\hat{s}}((\dot{a}(t)+(\dot{b}(t)-\dot{a}(t))\hat{s}\hat{x}(t,\hat{s}))/(b(t)-a(t)))$-term captures the {\bf underlying translation} of any point of interest $s(t,\hat{s})$ due to the moving boundaries;
\end{enumerate}

In the next sections, we use these dynamics to derive the power balance for the system~\eqref{eq:pH-DPS} on $\mathbb{S}_{ab}(t)$ and, subsequently, derive a corresponding Stokes-Dirac structure and pH formulation for the generic state evolution~\eqref{eq:dyn}. Before doing so, we briefly illustrate~\eqref{eq:dyn} for the telegrapher's equation.

\subsection{Illustrative example: Ideal transmission line}\label{sec:example}
    The dynamics of an ideal transmission line are governed by the telegrapher's equations
    \begin{equation*}
        \partial_t q(t,s) = -\partial_s I(t,s)\text{ and }\partial_t \phi(t,s) = -\partial_s V(t,s)
    \end{equation*}
    where $q(t,s)$ and $\phi(t,s)$ denote, respectively, the charge and flux density at time $t\in\mathbb{T}$ and position $s\in\mathbb{S}_{ab}(t)$. The current $I$ and voltage $V$ are given by 
    \begin{equation*}
        I(t,s) = \frac{\phi(t,s)}{L}\text{ and }V(t,s)=\frac{q(t,s)}{C},
    \end{equation*}
    with $L,C\in\mathbb{R}_{>0}$ being the distributed inductance and capacitance of the TL, which we assume to be constant. The state of the TL is $x=(q,\phi)\in\mathcal{X}(t)$, $t\in\mathbb{T}$, with the state space $\mathcal{X}~\colon~\mathbb{T}\rightarrow L_2([a(\cdot),b(\cdot)],\mathbb{R}^2)$. The stored energy in the segment $\mathbb{S}_{ab}(t)$ of the TL is given by the Hamiltonian $H~\colon~\mathcal{X}\times\mathbb{T}\rightarrow\mathbb{R}_{\geqslant 0}$ as 
    \begin{equation*}
        H(x,t) = \frac{1}{2}\int_{a(t)}^{b(t)} \frac{q^2(s)}{C} + \frac{\phi^2(s)}{L}\,\mathrm{d}s.
    \end{equation*}
    The effort variables of the TL are $e(t,s)=\delta_{x}H(x,t)=(I(t,s),V(t,s))$, and, hence, $Q=\operatorname{diag}(L^{-1},C^{-1})$.

    We apply the change of spatial variable~\eqref{eq:coc} and define the new state $\hat{x}\in\hat{\mathcal{X}}(t)$, $t\in\mathbb{T}$, with $\hat{\mathcal{X}}~\colon~\mathbb{T}\rightarrow L_2([0,1],\mathbb{R}^2)$, according to~\eqref{eq:xhat}, i.e., $\hat{x}(t,\hat{s})=(\hat{q}(t,\hat{s}),\hat{\phi}(t,\hat{s}))$ with $\hat{q}(t,\hat{s})\coloneqq \sqrt{b(t)-a(t)}q(t,h(t,\hat{s}))$ and $\hat{\phi}(t,\hat{s})\coloneqq \phi(t,h(t,\hat{s}))$. Using~\eqref{eq:dyn}, the dynamics of $\hat{x}$ are governed by
    \begin{align}\label{eq:TL-dyn}
        &\partial_t\begin{bmatrix}
            \hat{q}(t,\hat{s})\\
            \hat{\phi}(t,\hat{s})
        \end{bmatrix} = \begin{bmatrix}
            -\frac{1}{b(t)-a(t)}\partial_{\hat{s}}\hat{I}(t,\hat{s})\\
            -\frac{1}{b(t)-a(t)}\partial_{\hat{s}}\hat{V}(t,\hat{s})
        \end{bmatrix} +\\
        &\begin{bmatrix}
            -\frac{1}{2}\frac{\dot{b}(t)-\dot{a}(t)}{b(t)-a(t)}\hat{q}(t,\hat{s})+\partial_{\hat{s}}\left(\frac{\dot{a}(t)+(\dot{b}(t)-\dot{a}(t))\hat{s}}{b(t)-a(t)}\frac{\hat{V}(t,\hat{s})}{C}\right)\\
             -\frac{1}{2}\frac{\dot{b}(t)-\dot{a}(t)}{b(t)-a(t)}\hat{\phi}(t,\hat{s})+\partial_{\hat{s}}\left(\frac{\dot{a}(t)+(\dot{b}(t)-\dot{a}(t))\hat{s}}{b(t)-a(t)}\frac{\hat{I}(t,\hat{s})}{L}\right) 
        \end{bmatrix}.\nonumber
    \end{align} 
    where $\hat{e}(t,\hat{s}) = (\hat{I}(t,\hat{s}),\hat{V}(t,\hat{s}))=\delta_{\hat{x}}\hat{H}(\hat{x},t) = Q\hat{x}(t,\hat{s})$ with $\hat{H}$ as in~\eqref{eq:Hhat}. As discussed before, the time-derivatives of the charge and flux densities depend not only on the spatial derivatives of the current $I$ and voltage $V$, respectively, but also on the rate at which the charge and flux, respectively, in the interval $\mathbb{S}_{ab}(t)$ are (de)compressed and the rate at which the point of interest $s(t,\hat{s})$ in~\eqref{eq:coc} moves through the charge and flux field. We now show that the total charge and flux in $\mathbb{S}_{ab}(t)$ are conserved quantities. To illustrate this, we consider the total charge $\hat{Q}~\colon~\hat{\mathcal{X}}\times\mathbb{T}\rightarrow \mathbb{R}$ at time $t\in\mathbb{T}$ in $\mathbb{S}_{ab}(t)$, which, using~\eqref{eq:coc}, is given by
    \begin{equation*}
        \hat{Q}(\hat{x},t) = \int_0^1 \sqrt{b(t)-a(t)}\hat{q}(\hat{s})\,\mathrm{d}\hat{s}.
    \end{equation*}
    Differentiating $\hat{Q}(\hat{x},t)$ with respect to time yields 
    \begin{align*}
        &\frac{\mathrm{d}\hat{Q}(\hat{x}(t,\cdot),t)}{\mathrm{d}t} = \int_0^1 \frac{1}{2}\frac{\dot{b}(t)-\dot{a}(t)}{\sqrt{b(t)-a(t)}}\hat{q}(t,\hat{s}) +\\
        &\quad \sqrt{b(t)-a(t)}\partial_t\hat{q}(t,\hat{s})\,\mathrm{d}\hat{s}\\
        &= \int_0^1 \partial_{\hat{s}}\left(\frac{\dot{a}(t)+(\dot{b}(t)-\dot{a}(t))\hat{s}}{\sqrt{b(t)-a(t)}}\hat{q}(t,\hat{s})\right) -\\
        &\quad \frac{\partial_{\hat{s}}\hat{I}(t,\hat{s})}{\sqrt{b(t)-a(t)}}\,\mathrm{d}\hat{s},\\
        &= \left[\frac{\dot{a}(t)+(\dot{b}(t)-\dot{a}(t))\hat{s}}{\sqrt{b(t)-a(t)}}\hat{q}(t,\hat{s}) - \frac{\hat{\phi}(t,\hat{s})}{L\sqrt{b(t)-a(t)}}\right]_{\hat{s}=0}^{\hat{s}=1},
    \end{align*}
    which shows that charge only enters $\mathbb{S}_{ab}(t)$ at the boundaries. In fact, the movement of the boundaries effectively introduces an additional current equal to the charge density at the boundary multiplied with the velocity of the boundary. The same applies to the total flux in $\mathbb{S}_{ab}(t)$.

\subsection{Conservation of energy}\label{sec:energy}
Consider again the new energy variables $\hat{x}$ introduced in~\eqref{eq:xhat} in the previous section as well as the redefined Hamiltonian $\hat{H}$ in~\eqref{eq:H-coc} as a function of this new state $\hat{x}$. Clearly, $\hat{H}$ and $H$ are equivalent representations of the energy contained in the time-varying spatial domain $\mathbb{S}_{ab}(t)$. Hence, we can differentiate $\hat{H}(\hat{x}(t,\cdot))$ to obtain the power balance equation for the system~\eqref{eq:pH-DPS} on $\mathbb{S}_{ab}(t)$ as follows.
\setcounter{thm}{\thelemma}\stepcounter{lemma}
\begin{lem}\label{lem:power-balance}
    Consider the linear boundary pH system~\eqref{eq:pH-DPS} with spatial domain $\mathbb{S}_{ab}(t)$ and the Hamiltonian $\hat{H}$ in~\eqref{eq:Hhat}. Then, the corresponding power balance is given by
    \begin{align}
        &\int_0^1 \hat{e}^\top(t,\hat{s})\partial_t\hat{x}(t,\hat{s})\,\mathrm{d}\hat{s} = -\left[\frac{\hat{e}^\top(t,\hat{s})M^\top S_1M\hat{e}(t,\hat{s})}{b(t)-a(t)}\right]_{\hat{s}=0}^{\hat{s}=1} \nonumber\\
    &\,+\left[\frac{1}{2}\frac{\dot{a}(t)+(\dot{b}(t)-\dot{a}(t))\hat{s}}{b(t)-a(t)}\hat{e}^\top(t,\hat{s})\hat{x}(t,\hat{s})\right]_{\hat{s}=0}^{\hat{s}=1}.\label{eq:power-balance}
    \end{align}
    where $S_1=S_1^\top\in\mathbb{R}^{r\times r}$ and $M\in\mathbb{R}^{r\times n}$ are the matrices introduced in Lemma~\ref{lem:stokes-dirac}.
\end{lem}
\begin{proof}
    See Appendix~\ref{app:power-balance} for a detailed proof.
\end{proof}
The power balance equation in Lemma~\ref{lem:power-balance} captures the conservation of energy in the system. If $\dot{a}(t)=\dot{b}(t)=0$, we recover the standard power balance of the system~\eqref{eq:pH-DPS} (on a static spatial domain). As a result of the time-varying spatial domain $\mathbb{S}_{ab}(t)$, the power balance~\eqref{eq:power-balance} features an additional boundary term $\frac{1}{2}(\dot{b}(t)\hat{e}^\top(t,1)\hat{x}(t,1) - \dot{a}(t)\hat{e}^\top(t,0)\hat{x}(t,0)$. This additional boundary term captures the power flow into or out of the system as a result of the movement of the left and right boundary. 

Next, we illustrate these developments by revisiting the telegrapher's equation from Section~\ref{sec:example}.
\subsection{Illustrative example continued: Ideal TL}
    Consider again the ideal TL for which, in the previous section, we carried out the change of spatial variable~\eqref{eq:coc} and, subsequently, derived the dynamics of the new state $\hat{x}$ as well as the conservation of the total charge (and flux) in $\mathbb{S}_{ab}(t)$. We will now use Lemma~\ref{lem:power-balance} to investigate the conservation of energy in the segment $\mathbb{S}_{ab}(t)$ of the TL. To this end, we recall that $\hat{x}(t,\hat{s})=(\hat{q}(t,\hat{s}),\hat{\phi}(t,\hat{s}))$ and $\hat{e}(t,\hat{s})=(\hat{V}(t,\hat{s}),\hat{I}(t,\hat{s}))$. Moreover, since 
    \begin{equation*}
        J_1 = \begin{bmatrix}
            0 & -1\\
            -1 & 0
        \end{bmatrix}
    \end{equation*}
    is full rank, we have that $M=I_2$ and $S_1 = -\frac{1}{2}J_1$ as shown in~\cite{LeGorrec2005}. Then, substitution in~\eqref{eq:power-balance} yields the power balance
    \begin{align}
        &(b(t)-a(t))\int_0^1 \hat{V}(t,\hat{s})\partial_{t}\hat{q}(t,\hat{s}) + \hat{I}(t,\hat{s})\partial_t\hat{\phi}(t,\hat{s})\,\mathrm{d}\hat{s} = \nonumber\\
        &\quad\hat{V}(t,0)\hat{I}(t,0)+\dot{b}(t)\left(\frac{\hat{q}^2(t,1)}{2C}+\frac{\hat{\phi}^2(t,1)}{2L}\right)\label{eq:TL-power-balance}\\
        &\quad -\hat{V}(t,1)\hat{I}(t,1)- \dot{a}(t)\left(\frac{\hat{q}^2(t,0)}{2C}+\frac{\hat{\phi}^2(t,0)}{2L}\right).\nonumber
    \end{align}
    Note that, in addition to the familiar power exchange equal to the product of current and voltage at the boundaries, additional flows of power are introduced that are equal to the product of velocity of the boundary and the energy density at said boundary. Naturally, due to the lossless nature of the TL (and the lack of a distributed source), power only enters the time-varying segment $\mathbb{S}_{ab}(t)$ through its boundaries. 

\subsection{Port-Hamiltonian formulation}
In this section, we will present the main contribution of this paper, which consists of a port-Hamiltonian formulation of the boundary pH-DPS~\eqref{eq:pH-DPS} on a time-varying spatial domain $\mathbb{S}_{ab}(t)$ resulting in the novel class of \emph{moving-boundary port-Hamiltonian systems}. To this end, we first construct the associated efforts, flows, boundary ports, and the underlying Dirac structure.

\subsubsection{Flows and efforts}
The new co-energy variables, or efforts, $\hat{e}$ were already defined in~\eqref{eq:efforts}. Let $\hat{f}\in\hat{\mathcal{X}}(t)$, $t\in\mathbb{T}$, denote the corresponding flows given by
\begin{align}
    &\hat{f}(t,\hat{s})\coloneqq \partial_t\hat{x}(t,\hat{s}) \stackrel{\eqref{eq:dyn}}{=}\hat{\mathcal{J}}(t)\hat{e}(t,\hat{s}) -\label{eq:flows}\\
    &\frac{1}{2}\frac{\dot{b}(t)-\dot{a}(t)}{b(t)-a(t)}Q^{-1}\hat{e}(t,\hat{s})+\partial_{\hat{s}}\left(\frac{\partial_th(t,\hat{s})}{b(t)-a(t)}Q^{-1}\hat{e}(t,\hat{s})\right).\nonumber
\end{align}

\subsubsection{Boundary ports}
Considering the power balance in~\eqref{eq:power-balance}, we introduce the following boundary ports $\hat{e}_\partial,\hat{f}_\partial~\colon~\mathbb{T}\rightarrow \mathbb{R}^{r}\times \mathbb{K}^n$ given by
\begin{equation}
    \begin{bmatrix}
        \hat{f}_{\partial}(t)\\
        \hat{e}_{\partial}(t)
    \end{bmatrix} = \begin{bmatrix}
        S_1M & -S_1M\\
        -\overline{\sqrt{\dot{a}(t)}}I_n & \overline{\sqrt{\dot{b}(t)}}I_n\\
        M & M\\
        \frac{1}{2}\sqrt{\dot{a}(t)}Q^{-1} & \frac{1}{2}\sqrt{\dot{b}(t)}Q^{-1}\\
    \end{bmatrix}\frac{(\hat{e}(t,0),\hat{e}(t,1))}{\sqrt{b(t)-a(t)}},
    \label{eq:boundary-ports}
\end{equation}
where $\overline{z}$ denotes the complex conjugate of $z\in\mathbb{C}$, such that~\eqref{eq:power-balance} reads $\langle \hat{f},\hat{e}\rangle_{L_2} = \langle \hat{f}_{\partial},\hat{e}_{\partial}\rangle_{\mathbb{R}^{r}\times\mathbb{K}^{n}}$. While the individual port variables become imaginary when $\dot{a}(t)<0$ or $\dot{b}(t)<0$, their product $\hat{f}_{\partial}^\hop(t)\hat{e}_{\partial}(t)$ always remains real-valued and, thus, preserves the classical interpretation as physical power flowing into or out of the spatial domain $\mathbb{S}_{ab}(t)$. 

 \subsubsection{Dirac structure}
 We note that in~\eqref{eq:flows} $\hat{x}=Q^{-1}\hat{e}$ and, similarly, in~\eqref{eq:boundary-ports} $\tr\hat{x}=Q^{-1}\tr\hat{e}$. We express these in terms of $\hat{e}$ to facilitate the construction of a new \emph{time-varying} Stokes-Dirac structure for the boundary pH-DPS~\eqref{eq:pH-DPS} on $\mathbb{S}_{ab}(t)$ as follows.
\setcounter{thm}{\thetheorem}\stepcounter{theorem}
\begin{thm}\label{thm:novel-boundary-stokes-dirac}
    Let $r=\rank J_1$, $\hat{\mathcal{F}}~\colon~\mathbb{T}\rightarrow L_2([0,1],\mathbb{R}^{n})\times \mathbb{R}^{r}\times \mathbb{K}^n$, and $\hat{\mathcal{E}}(t)=\hat{\mathcal{F}}^*(t)$ for all $t\in\mathbb{T}$. Let $\hat{\mathcal{B}}(t)=\hat{\mathcal{F}}(t)\times\hat{\mathcal{E}}(t)$ be equipped with the non-degenerate bilinear form $\left\langle\cdot\middle|\cdot\right\rangle_{\hat{\mathcal{B}}}~\colon~\hat{\mathcal{B}}(t)\rightarrow\mathbb{K}$ given by
    \begin{equation*}
        \langle (\hat{f},\hat{f}_{\partial})\mid(\hat{e},\hat{e}_{\partial})\rangle_{\hat{\mathcal{B}}} = \langle \hat{f},\hat{e}\rangle_{L_2} - \langle \hat{f}_{\partial},\hat{e}_{\partial}\rangle_{\mathbb{R}^{r}\times\mathbb{K}^{n}}.
    \end{equation*}
    Then, there exist a symmetric matrix $S_1\in\mathbb{R}^{r\times r}$ and a full rank matrix $M\in\mathbb{R}^{r\times n}$ such that the vector subspace 
    \begin{align*}
        &\hat{\mathcal{D}}(t)\coloneqq \left\{(\hat{f},\hat{f}_{\partial},\hat{e},\hat{e}_{\partial}\in\hat{\mathcal{B}}(t)\middle|\right.\hat{e}\in\hat{\mathcal{E}}(t),\partial_{\hat{s}}\hat{e}\in\hat{\mathcal{F}}(t),\\
        &~\hat{f}(t,\hat{s})=-\frac{1}{2}\frac{\dot{b}(t)-\dot{a}(t)}{b(t)-a(t)}Q^{-1}\hat{e}(t,\hat{s}) +\\
        &\quad \partial_{\hat{s}}\left(\frac{\dot{a}(t)+(\dot{b}(t)-\dot{a}(t))\hat{s}}{b(t)-a(t)}Q^{-1}\hat{e}(t,\hat{s})\right) + \mathcal{\hat{J}}(t)\hat{e}(t,\hat{s}),\\
        &\left.\begin{bmatrix}
            \hat{f}_{\partial}(t)\\
            \hat{e}_{\partial}(t)
        \end{bmatrix} = \begin{bmatrix}
            S_1M & -S_1M\\
            -\overline{\sqrt{\dot{a}(t)}}I_n & \overline{\sqrt{\dot{b}(t)}}I_n\\
            M & M\\
            \frac{1}{2}\sqrt{\dot{a}(t)}Q^{-1} & \frac{1}{2}\sqrt{\dot{b}(t)}Q^{-1}
        \end{bmatrix}\frac{(\hat{e}(t,0),\hat{e}(t,1))}{\sqrt{b(t)-a(t)}}\right\},
    \end{align*}
    is, for all $t\in\mathbb{T}$, a Dirac structure with respect to the symmetric bilinear form $\left\llangle\cdot,\cdot\right\rrangle_{\hat{\mathcal{B}}}~\colon~\hat{\mathcal{B}}(t)\times\hat{\mathcal{B}}(t)\rightarrow\mathbb{K}$ given by
    \begin{align*}
        &\left\llangle \left(\hat{f}^1,\hat{f}^1_{\partial},\hat{e}^1,\hat{e}^1_{\partial}\right),\left(\hat{f}^2,\hat{f}^2_{\partial},\hat{e}^2,\hat{e}^2_{\partial}\right)\right\rrangle_{\hat{\mathcal{B}}} = \\
        &\quad\left\langle \left(\hat{f}^1,\hat{f}^1_{\partial}\right)\middle|\left(\hat{e}^2,\hat{e}^2_{\partial}\right)\right\rangle_{\hat{\mathcal{B}}} +\left\langle \left(\hat{f}^2,\hat{f}^2_{\partial}\right)\middle|\left(\hat{e}^1,\hat{e}^1_{\partial}\right)\right\rangle_{\hat{\mathcal{B}}}.
    \end{align*}
\end{thm}
\begin{proof}
    See Appendix~\ref{app:novel-boundary-stokes-dirac} for a detailed proof.
\end{proof}
The Dirac structure in Theorem~\ref{thm:novel-boundary-stokes-dirac} is essentially a combination of Stokes theorem, which describes the internal dynamics of the boundary pH-DPS (see, e.g.,~\cite{LeGorrec2005}), and Leibniz integral rule (see, e.g.,~\cite{Protter1985}), which captures the effect of the moving boundaries, formulated in terms of the boundary ports introduced earlier. 

We note that Theorem~\ref{thm:novel-boundary-stokes-dirac} exploits the additional generality of Definition~\ref{def:Dirac} to incorporate the complex-valued boundary ports related to the movement of the boundaries introduced in~\eqref{eq:boundary-ports}. 

Using the time-varying Dirac structure $\hat{\mathcal{D}}(t)$ formulated in Theorem~\ref{thm:novel-boundary-stokes-dirac}, we can now formally establish the class of \emph{moving-boundary port-Hamiltonian systems}, as follows.
\setcounter{thm}{\thedefinition}\stepcounter{definition}
\begin{prop}\label{prop:mbphs}
    The moving-boundary pH system on the state space $\hat{\mathcal{X}}~\colon~\mathbb{T}\rightarrow L_2([0,1],\mathbb{R}^{n})$ with respect to the Dirac structure $\hat{\mathcal{D}}(t)$ associated with $\hat{\mathcal{J}}(t)$ and with $a(t)$ and $b(t)$ satisfying Assumption~\ref{asm:bounds}, and generated by the Hamiltonian functional $\hat{H}$ in~\eqref{eq:Hhat}, is the dynamical system 
\begin{equation}
    (\partial_t\hat{x},\hat{f}_{\partial},Q\hat{x},\hat{e}_{\partial})\in\hat{\mathcal{D}}(t),\quad t\in\mathbb{T}.
\end{equation}
\end{prop}
Observe that the above definition of a moving-boundary pH system genuinely generalizes Definition~\ref{def:lin-fo-bound-pH-DPSs}, i.e., Proposition~\ref{prop:mbphs} reduces to Definition~\ref{def:lin-fo-bound-pH-DPSs} when $\dot{a}(t)=\dot{b}(t)=0$.

\subsection{Illustrative example continued: Ideal TL}
    Consider once more the ideal TL from Section~\ref{sec:example}, for which the boundary ports~\eqref{eq:boundary-ports} are
    \begin{align*}
        &\left[\begin{array}{@{}c;{2pt/2pt}c@{}}\
            \hat{f}_{\partial}(t) & \hat{e}_{\partial}(t)
        \end{array}\right] = \frac{1}{\sqrt{b(t)-a(t)}}\times\\
        &\left[\begin{array}{@{}c;{2pt/2pt}c@{}}
            \frac{1}{2}(\hat{I}(t,0)-\hat{I}(t,1)) & \hat{V}(t,0)+\hat{V}(t,1)\\
            \frac{1}{2}(\hat{V}(t,0)-\hat{V}(t,1)) & \hat{I}(t,0) + \hat{I}(t,1)\\
            \overline{\sqrt{\dot{b}(t)}}\hat{V}(t,1)-\overline{\sqrt{\dot{a}(t)}}\hat{V}(t,0) & \frac{\sqrt{\dot{a}(t)}\hat{q}(t,0)+\sqrt{\dot{b}(t)}\hat{q}(t,1)}{2}\\
            \overline{\sqrt{\dot{b}(t)}}\hat{I}(t,1)-\overline{\sqrt{\dot{a}(t)}}\hat{I}(t,0) & \frac{\sqrt{\dot{a}(t)}\hat{\phi}(t,0)+\sqrt{\dot{b}(t)}\hat{\phi}(t,1)}{2}
        \end{array}\right].
    \end{align*}
    The dynamics of the segment $\mathbb{S}_{ab}(t)$ can be represented as the \emph{moving-boundary port-Hamiltonian system} $(\partial_t\hat{x},\hat{f}_{\partial},Q\hat{x},\hat{e}_{\partial})\in\hat{\mathcal{D}}(t)$ with $\hat{\mathcal{D}}(t)$ as in Theorem~\ref{thm:novel-boundary-stokes-dirac}.

\section{Application to dynamic meshing}\label{sec:numerical-results}
\subsection{Spatial discretization scheme}\label{sec:discretization-scheme}
In this section, we use Theorem~\ref{thm:novel-boundary-stokes-dirac} to discretize a moving-boundary pH system using a dynamic grid. To this end, we revisit the lossless TL considered in Section~\ref{sec:example} and exploit our results in Theorem~\ref{thm:novel-boundary-stokes-dirac} to modify the discretization scheme proposed in~\cite{Golo2004}. The resulting discretization scheme is a simple mixed-element method that is capable of incorporating the time-varying nature of the spatial domain by adopting a dynamic mesh. 

It is important to note that this modification means that the discretization is no longer structure-preserving, however, our pH formulation and the proposed discretization scheme recover the results of~\cite{Golo2004} when, for some period of time, the spatial domain is time-invariant.

Consider the time-varying segment $\mathbb{S}_{ab}(t)=[a(t),b(t)]$ of the lossless TL from Section~\ref{sec:example}. We discretize this system by subdividing $\mathbb{S}_{ab}(t)$ into $N$ time-varying elements with the $N+1$ dynamic grid points $\{s_i(t)\}_{i\in\mathcal{N}\cup (N+1)}$, $\mathcal{N}\coloneqq \{1,2,\hdots,N\}$, with \begin{equation}\label{eq:dyn-grid} s_i = a(t) + \frac{i(b(t)-a(t))}{N},\quad \text{ for }i\in\mathcal{N}\cup (N+1),\end{equation} being the boundary between the $(i-1)$-th and $i$-th element. Through the time-varying change of coordinates~\eqref{eq:coc} this dynamic gridding of $\mathbb{S}_{ab}(t)$ maps to a static gridding of the transformed spatial domain $\hat{s}\in[0,1]$. Precisely, the dynamic grid points $\{s_i(t)\}_{i\in\mathcal{N}\cup (N+1)}$ in~\eqref{eq:dyn-grid} correspond to the static grid points $\{\hat{s}_i\}_{i\in\mathcal{N}\cup(N+1)}$ given by \[\hat{s}_i=\frac{(i-1)}{N},\quad\text{for }i\in\mathcal{N}\cup(N+1).\] These static grid points subdivide the transformed spatial domain $[0,1]$ into $N$ constant elements of equal length.

Our discretization scheme is based on approximating the specific Dirac structure in Theorem~\ref{thm:novel-boundary-stokes-dirac} using this static grid of the transformed spatial domain $\hat{s}\in [0,1]$. Following the same procedure as in~\cite{Golo2004}, we approximate the energy variables $\hat{x}(t,\hat{s})=(\hat{q}(t,\hat{s}),\hat{\phi}(t,\hat{s}))$ according to
\begin{equation*}
    \hat{x}(t,\hat{s}) \approx \sum_{\mathclap{i\in\mathcal{N}}} \hat{x}_{i,i+1}(t)\omega_{i,i+1}(\hat{s}),
\end{equation*}
where, for all $i\in\mathcal{N}$,
\begin{equation*}
    \omega_{i,i+1}(\hat{s}) = \begin{cases}
        N,\qquad &\text{for }\hat{s}\in(\hat{s}_i,\hat{s}_{i+1}),\\
        N/2,\qquad &\text{for }s\in\{\hat{s}_i,\hat{s}_{i+1}\},\\
        0,\qquad &\text{elsewhere.}
    \end{cases}
\end{equation*}
Similarly, we approximate the flows $\hat{f}(t,\hat{s}) = \partial_{t}\hat{x}(t,\hat{s})$ as 
\begin{equation*}
    \hat{f}(t,\hat{s}) \approx \sum_{\mathclap{i\in\mathcal{N}\cup\{N+1\}}}\hat{f}_{i,i+1}(t)\omega_{i,i+1}(\hat{s}).
\end{equation*}
The effort variables $\hat{e}=(\hat{V},\hat{I})$, on the other hand, are approximated as
\begin{equation}
    \hat{e}(t,\hat{s}) \approx \sum_{\mathclap{i\in\mathcal{N}\cup\{N+1\}}} \hat{e}_i(t)\omega_{i}(\hat{s}),
\end{equation}
where, for all $i\in\mathcal{N}\cup\{N+1\}$,
\begin{equation*}
    \omega_i(\hat{s}) = \begin{cases}
        N(\hat{s}_{i+1}-\hat{s}),\quad &s\in[\hat{s}_i,\hat{s}_{i+1}),\\
        N(\hat{s}-\hat{s}_{i-1}),\quad &s\in(\hat{s}_{i-1},\hat{s}_i),\\
        0,\quad &\text{otherwise.}
    \end{cases}
\end{equation*}
Note that these basis functions satisfy
\begin{equation*}
    \omega_{i,i+1}(\hat{s}) = -\partial_{\hat{s}}\omega_i(\hat{s}) = \partial_{\hat{s}}\omega_{i+1}(\hat{s})\text{ for }\hat{s}\in(\hat{s}_i,\hat{s}_{i+1}),
\end{equation*}
for all $i,\in\mathcal{N}$, and, for all $i,j\in\mathcal{N}\cup\{N+1\}$,
\begin{equation*}
    \omega_i(\hat{s}_j) = \begin{cases}
        1,\quad\text{if }i=j,\\
        0,\quad\text{if }i\neq j.
    \end{cases}
\end{equation*}

Next, we will use these approximations to approximate the Dirac structure $\hat{\mathcal{D}}(t)$ for the lossless TL. To this end, we substitute the approximate energy variables, efforts and flows into the TL dynamics~\eqref{eq:TL-dyn} and integrate over $[\hat{s}_i,\hat{s}_{i+1}]$ to obtain, for all $i\in\mathcal{N}$,
\begin{align*}
    0 &= \int_{\hat{s}_i}^{\hat{s}_{i+1}}\hat{f}_{i,i+1}(t)\omega_{i,i+1}(\hat{s}) +\\
    &\frac{1}{2}\frac{\dot{b}(t)-\dot{a}(t)}{b(t)-a(t)}\hat{x}_{i,i+1}(t)\omega_{i,i+1}(\hat{s})\,\mathrm{d}\hat{s} +\\
    &\sum_{i=1}^{N+1}\left[\frac{1}{b(t)-a(t)}\begin{bmatrix} 0 & 1\\
    1 & 0\end{bmatrix}\omega_i(\hat{s})\right]_{\hat{s}_i}^{\hat{s}_{i+1}}\hat{e}_i(t)-\\
    &\sum_{i=1}^{N+1}\left[\frac{\dot{a}(t)+(\dot{b}(t)-\dot{a}(t))\hat{s}}{b(t)-a(t)}\omega_i(\hat{s})\right]_{\hspace*{.1cm}\mathclap{\hat{s}_{i}}}^{\hspace*{.2cm}\mathclap{\hat{s}_{i+1}}}\hat{e}_i(t)=\\
    &\hat{f}_{i}(t) + \frac{1}{2}\frac{\dot{b}(t)-\dot{a}(t)}{b(t)-a(t)}\hat{x}_{i,i+1}(t) + \\
    &\begin{bmatrix}
        \frac{\dot{a}(t)+(\dot{b}(t)-\dot{a}(t))\hat{s}_{i}}{b(t)-a(t)} & -\frac{1}{b(t)-a(t)}\\
        -\frac{1}{b(t)-a(t)} & \frac{\dot{a}(t)+(\dot{b}(t)-\dot{a}(t))\hat{s}_{i}}{b(t)-a(t)}
    \end{bmatrix}\hat{e}_{i}(t)-\\
    &\begin{bmatrix}
        \frac{\dot{a}(t)+(\dot{b}(t)-\dot{a}(t))\hat{s}_{i+1}}{b(t)-a(t)} & -\frac{1}{b(t)-a(t)}\\
        -\frac{1}{b(t)-a(t)} & \frac{\dot{a}(t)+(\dot{b}(t)-\dot{a}(t))\hat{s}_{i+1}}{b(t)-a(t)}
    \end{bmatrix}\hat{e}_{i+1}(t),
\end{align*}
where, for all $i\in\mathcal{N}$,
\begin{equation*}
    \dot{\hat{x}}_{i,i+1}(t) = \hat{f}_{i,i+1}(t).
\end{equation*}
Similarly, we substitute the approximated efforts and flows into the power balance~\eqref{eq:TL-power-balance} 
\begin{align*}
    0 &= \sum_{i\in\mathcal{N}} \frac{1}{2}\frac{\dot{a}(t)\hat{e}_i^\top(t)Q^{-1}\hat{e}_i(t)-\dot{b}(t)\hat{e}_{i+1}^\top(t) Q^{-1}\hat{e}_{i+1}(t)}{b(t)-a(t)}+\\
    &\frac{1}{b(t)-a(t)}\hat{e}_{i+1}^\top(t)\begin{bmatrix} 0 & I\\
    I & 0\end{bmatrix}\hat{e}_{i+1}(t) -\\
    &\frac{1}{b(t)-a(t)}\hat{e}_{i}^\top(t)\begin{bmatrix} 0 & I\\
    I & 0\end{bmatrix}\hat{e}_{i}(t)+\\
    &\hat{f}_{i,i+1}^\top(t)\int_{\hat{s}_i}^{\hat{s}_{i+1}} (\hat{e}_i(t)\omega_i(\hat{s}) + \hat{e}_{i+1}(t)\omega_{i+1}(\hat{s}))\omega_{i,i+1}(\hat{s})\,\mathrm{d}\hat{s}\\
    &= \sum_{i\in\mathcal{N}} \hat{f}_{i,i+1}^\top(t)\hat{e}_{i,i+1}(t) - (\hat{f}^{\partial}_{i,i+1}(t))^\top\hat{e}_{i,i+1}^\partial(t),
    \end{align*}
where we have introduced the effort corresponding to $\hat{f}_{i,i+1}(t)$ as
\begin{align*}
    \hat{e}_{i,i+1}(t) &\coloneqq \int_{\hat{s}_i}^{\hspace*{.2cm}\mathclap{\hat{s}_{i+1}}} (\hat{e}_i(t)\omega_i(\hat{s}) + \hat{e}_{i+1}(t)\omega_{i+1}(\hat{s}))\omega_{i,i+1}(\hat{s})\,\mathrm{d}\hat{s}\\
    &= \frac{1}{2}(\hat{e}_i(t)+\hat{e}_{i+1}(t)),
\end{align*}
for all $i\in\mathcal{N}$, and the boundary ports, for $i\in\mathcal{N}$,
\begin{align*}
    &\hat{f}^\partial_{i,i+1}(t) \coloneqq \frac{1}{\sqrt{b(t)-a(t)}}\times\\
    &\begin{bmatrix}\begin{smallmatrix}
        \frac{1}{2}\begin{bmatrix}\begin{smallmatrix}
            0 & 1\\
            1 & 0
        \end{smallmatrix}\end{bmatrix}(\hat{e}_{i}(t)-\hat{e}_{i+1}(t))\\
        \overline{\sqrt{\partial_th(t,\hat{s}_{i+1})}}\hat{e}_{i+1}(t)-\overline{\sqrt{\partial_th(t,\hat{s}_i)}}\hat{e}_i(t)
    \end{smallmatrix}\end{bmatrix},\\
    &\hat{e}^\partial_{i,i+1}(t) \coloneqq \frac{1}{\sqrt{b(t)-a(t)}}\times\\
    &\begin{bmatrix}\begin{smallmatrix}
        \hat{e}_{i}(t)+\hat{e}_{i+1}(t)\\
        \frac{1}{2}Q^{-1}(\sqrt{\partial_th(t,\hat{s}_i)}\hat{e}_i(t) + \sqrt{\partial_th(t,\hat{s}_{i+1})}\hat{e}_{i+1}(t))
    \end{smallmatrix}\end{bmatrix}.
\end{align*}
Note that the discretization of the power balance differs from~\cite{Golo2004} to incorporate the moving boundaries and, as a result, the discretization scheme does not preserve the power balance of the moving-boundary pH system. However, when $\dot{a}(t)=\dot{b}(t)=0$, the discretization scheme proposed here recovers the one in~\cite{Golo2004}, and, as such, the scheme is structure-preserving in this specific case. 

Finally, the Hamiltonian is discretized as follows.
\begin{align*}
    &\hat{H}(\hat{x},t) \\
    &\quad=\frac{1}{2}\sum_{i\in\mathcal{N}} \int_{\hat{s}_i}^{\hat{s}_{i+1}} \hat{x}_{i,i+1}^\top(t)Q^{-1}\hat{x}_{i,i+1}(t)\omega_{i,i+1}^2(\hat{s})\,\mathrm{d}\hat{s},\\
    &\quad=\frac{1}{2}\sum_{i\in\mathcal{N}} \underbrace{N\hat{x}^\top_{i,i+1}(t) Q^{-1}\hat{x}_{i,i+1}(t)}_{\eqqcolon \hat{H}_i(\hat{x}_{i,i+1}(t))},
\end{align*}
which we use to obtain, for all $i\in\mathcal{N}$,
\begin{equation*}
    \hat{e}_{i,i+1}(t) = \partial_{\hat{x}_{i,i+1}} \hat{H}_{i}(\hat{x}_{i,i+1}) = NQ^{-1}\hat{x}_{i,i+1}(t).
\end{equation*}
At this point, the only remaining step to complete the discretization is to set the boundary conditions of the TL segment $\mathbb{S}_{ab}(t)$, which we will perform in the next section.

\subsection{Simulation results}
We now use the above discretization scheme to simulate a lossless TL on a time-varying spatial domain. We consider a TL with $L = 1$ \si{\henry\per\meter} and $C = 1$ \si{\farad\per\meter}. The spatial domain of the TL as a whole is $[0,1]$ \si{\meter} with the voltage at $s=0$ \si{\meter} being a sinusoidal input, i.e., $e^q(t,0) = u(t)=\sin(t)$ \si{\volt}, while the other end of the transmission line is truncated by a $R=1$ \si{\ohm} resistor such that $e^{\phi}(t,1)=e^q(t,1)$ for all $t\in\mathbb{R}$. The initial condition is $q(0,s)=\phi(0,s)=0$ for all $s\in[0,1]$. The analytic solution for the voltage and current are given by $V_{\mathrm{true}}(t,s) = I_{\mathrm{true}}(t,s) = \sin(\max\{0,t-s\})$ for all $t\in\mathbb{R}_{\geqslant 0}$ and $s\in[0,1]$.

\begin{figure}[!bt]
        \centering
        \begin{subfigure}{\linewidth}
            \includegraphics[width=\linewidth]{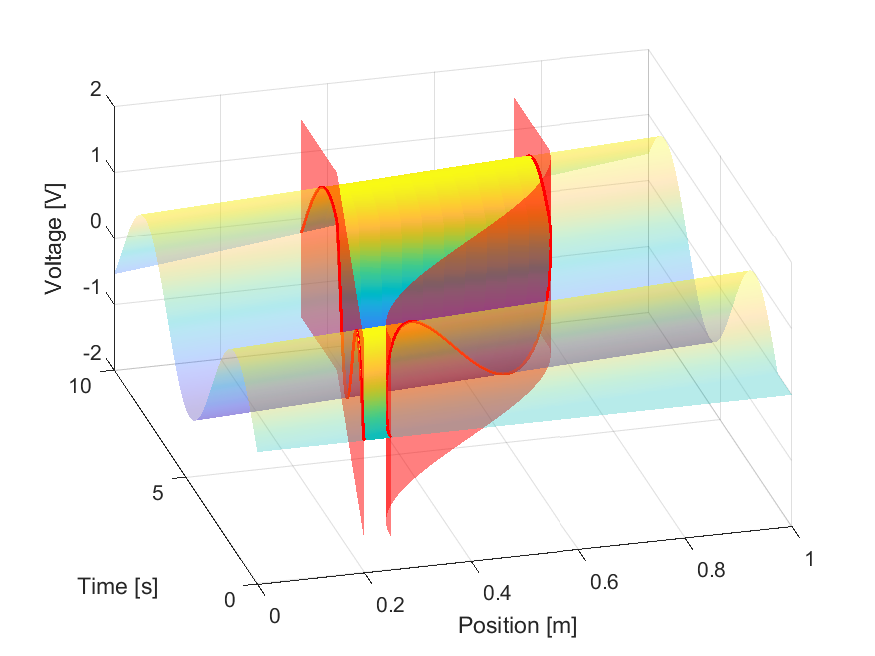}
            \caption{Side view.}
            \label{fig:simulation-side}
        \end{subfigure}        
        \begin{subfigure}{\linewidth}
            \includegraphics[width=\linewidth]{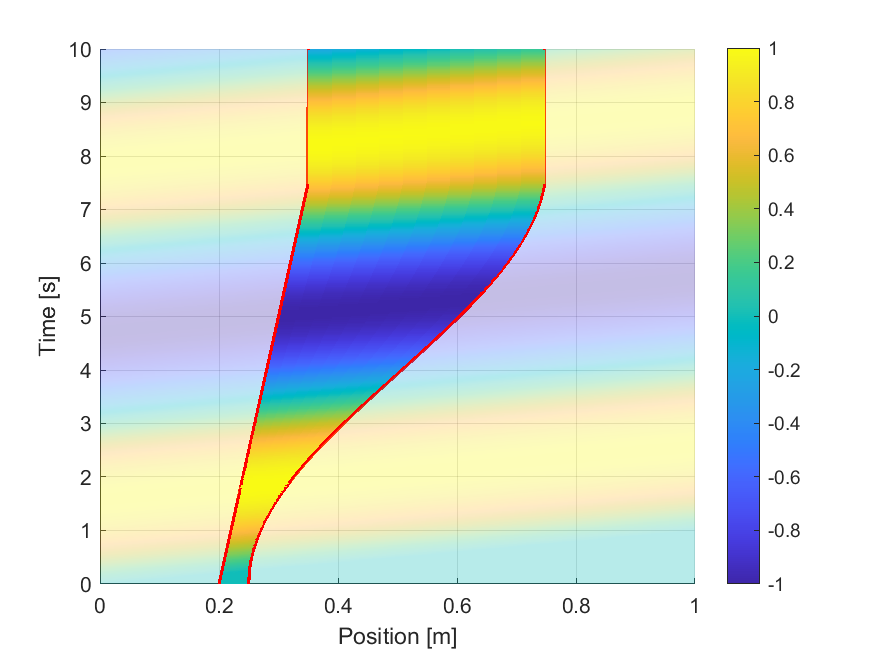}
            \caption{Top view.}
            \label{fig:simulation-top}
        \end{subfigure}        
        \caption{Simulated voltage (\protect\tikz\protect\shade[shading=parula,shading angle=135,opacity=1] (0,0) rectangle (0.2,0.2);) in the time-varying segment $\mathbb{S}_{ab}(t)$ (\protect\tikz\protect\filldraw[fill=red,opacity=0.5,draw=none] (0,0) rectangle (0.2,0.2);). Outside $\mathbb{S}_{ab}(t)$, the analytic solution $V_{\mathrm{true}}$ (\protect\tikz\protect\shade[shading=parula,shading angle=135,opacity=0.3] (0,0) rectangle (0.2,0.2);) is depicted.}
        \label{fig:simulation}
    \end{figure}
    \begin{figure}[!ht]
        \centering
        \begin{subfigure}{\linewidth}
            \includegraphics[width=\linewidth]{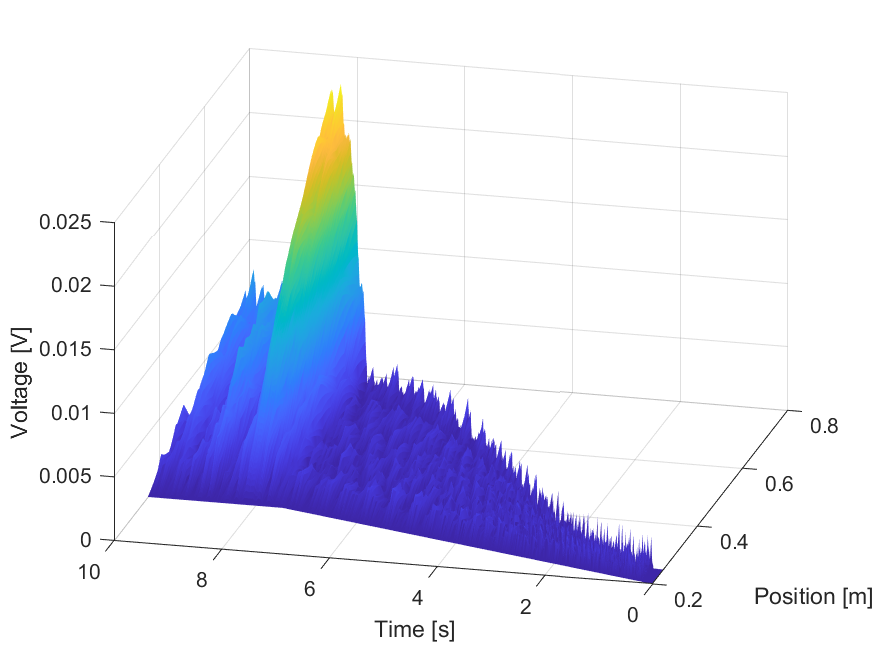}
            \caption{Side view.}
            \label{fig:error-side}
        \end{subfigure}        
        \begin{subfigure}{\linewidth}
            \includegraphics[width=\linewidth]{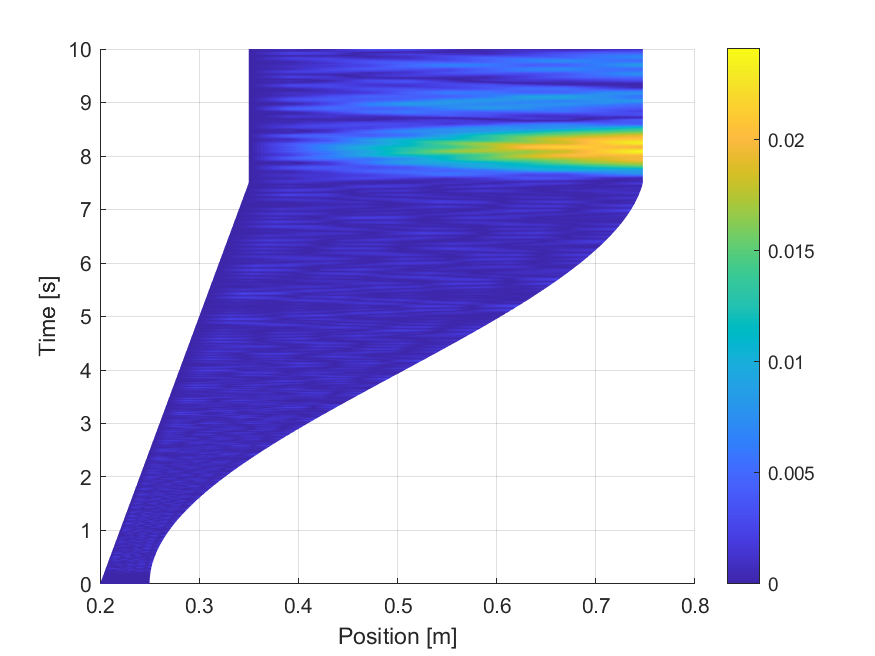}
            \caption{Top view.}
            \label{fig:error-top}
        \end{subfigure}     
        \caption{Absolute value of the error with respect to the analytic solution $V_{\mathrm{true}}$ (\protect\tikz\protect\shade[shading=parula,shading angle=135,opacity=1] (0,0) rectangle (0.2,0.2);) in the time-varying segment $\mathbb{S}_{ab}(t)$.}
        \label{fig:error}
    \end{figure}
    \begin{figure}[!ht]
        \centering
        \includegraphics[width=\linewidth]{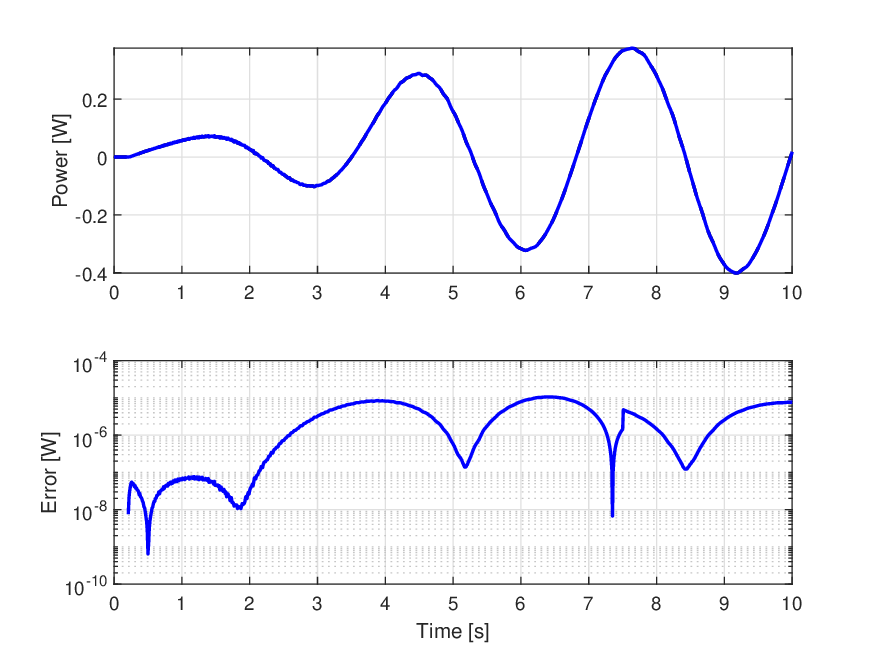}
        \caption{Power balance showing the power flowing into the TL segment $\mathbb{S}_{ab}(t)$ (\kern0.5pt\protect\tikz[scale=0.6]{\protect\draw[dashed, line width=1pt] (0,0.1)--(0.55,0.1);\protect\draw[color=white] (0,0)--(0.55,0)}) and the time-derivative of the discretized Hamiltonian (\kern0.5pt\protect\tikz[scale=0.6, line width=1pt]{\protect\draw[color=blue] (0,0.1)--(0.55,0.1);\protect\draw[color=white] (0,0)--(0.55,0)}), i.e., $\sum_{i=1}^{N} (\mathrm{d}\hat{H}_i(\hat{x}_{i,i+1}(t))/(\mathrm{d}t)$. The absolute error between the power flowing in and the change of the energy stored in the Hamiltonian is shown in the bottom plot.}
        \label{fig:power}
    \end{figure}

To illustrate our results, we consider the time-varying segment given by
\begin{equation*}
    \mathbb{S}_{ab}(t) = \begin{cases}
        [2+0.2t,5-\cos(0.25t)]/10,~\text{ for }t\leqslant \bar{t},\\
        [2+0.2\bar{t},5-\cos(0.25\bar{t})]/10,~\text{ for }t> \bar{t},\\
    \end{cases}
\end{equation*}
where $\bar{t}=7.5$ \si{\second}. In other words, $a(t) = 0.2+0.02t$ and $b(t) = 0.5-0.25\cos(0.25t)$ for $t\in[0,7.5]$ after which the boundaries remain static at $a(7.5)$ and $b(7.5)$. We use the discretization scheme proposed in Section~\ref{sec:discretization-scheme} to simulate this part of the TL using $N=10$ mixed elements. To impose boundary conditions, we partition the effort according to $\hat{e}(t,\hat{s})=(\hat{e}^{\hat{q}}(t,\hat{s}),\hat{e}^{\hat{\phi}}(t,\hat{s}))$ and we impose that the voltage at the boundary $\hat{s}=0$ (i.e., $s(t,0)=a(t)$) satisfies $\hat{e}^{\hat{q}}(t,0) = \sqrt{b(t)-a(t)}V_{\mathrm{true}}(t,a(t))$ for all $t\in\mathbb{R}_{\geqslant 0}$. Similarly, the current at the boundary $\hat{s}=1$  (i.e., $s(t,1)=b(t)$) is set to $\hat{e}^{\hat{\phi}}(t,1) = \sqrt{b(t)-a(t)}I_{\mathrm{true}}(t,b(t))$ for all $t\in\mathbb{R}_{\geqslant 0}$. 
    
The simulated voltage in $\mathbb{S}_{ab}(t)$ is depicted in Fig.~\ref{fig:simulation} along with the time-varying boundaries of $\mathbb{S}_{ab}(t)$. Outside $\mathbb{S}_{ab}(t)$, the analytic solution is shown as well. The simulation error is shown in Fig.~\ref{fig:error}, where we see that simulation closely resembles the analytic solution and is reasonably accurate. In fact, the error can be further reduced by using more than $N=10$ mixed elements. Fig.~\ref{fig:power} shows the time-derivative of the Hamiltonian, which corresponds to the instantaneous power being stored in the discretized TL segment $\mathbb{S}_{ab}(t)$, and the power entering through the boundary ports. We observe that, as already mentioned in the previous section, the modified discretization scheme is not structure-preserving, which manifests in the time-derivative of the discretized Hamiltonian not being exactly equal to the power entering through the boundary ports. In other words, the power balance of the underlying distributed parameter system is not preserved. However, when $\dot{a}(t)=\dot{b}(t)=0$, we see that the discretization scheme is structure-preserving (in fact, it is then equivalent to the scheme presented in~\cite{Golo2004}) such that increase or decrease of the Hamiltonian can be fully explained by the boundary ports.

\section{Conclusions}\label{sec:conclusions}
In this paper, we formally introduced the class of moving-boundary port-Hamiltonian systems. We showed that boundary pH-DPSs on a time-varying one-dimensional spatial domain give rise to a specific \emph{time-varying} Stokes-Dirac structure and admit a corresponding pH formulation. We studied the dynamics of this new class of systems as well as their conservation of energy. Finally, we demonstrated that our results can be leveraged to develop a spatial discretization scheme with dynamic meshing for approximating the telegrapher's equations on a time-varying spatial domain, which we subsequently verified numerically.

\bibliographystyle{plain}        
\bibliography{refs}           

\section*{Appendix}
\appendix
\section{Proof of Lemma~\ref{lem:dyn}}\label{app:dynamics}
We will now derive the dynamics that govern $\hat{x}\in\hat{\mathcal{X}}(t)$, $t\in\mathbb{T}$. To this end, we note that
\begin{equation}\label{eq:dsxhat}
    \partial_{\hat{s}}\hat{x}(t, \hat{s})=\left(b(t)-a(t)\right)^{\frac{3}{2}}\left(\partial_{\sigma}x(t,\sigma)\right)|_{\sigma=h(t,\hat{s})}.
\end{equation}
Using the chain rule, the dynamics of $\hat{x}$ are governed by
\begin{align}
    &\partial_t\hat{x}(t,\hat{s}) = \frac{1}{2}\frac{\dot{b}(t)-\dot{a}(t)}{\sqrt{b(t)-a(t)}}x(t,h(t,\hat{s})) + \sqrt{b(t)-a(t)}\times\nonumber\\
    &\quad \left((\partial_\sigma x(t,\sigma)\mid_{\sigma=h(t,\hat{s})}\partial_th(t,\hat{s}) + (\partial_tx(t,\sigma))\mid_{\sigma=h(t,\hat{s})}\right),\nonumber\\
    &\stackrel{\eqref{eq:pH-DPS}}{=} \left(\frac{1}{2}\frac{\dot{b}(t)-\dot{a}(t)}{\sqrt{b(t)-a(t)}} + \sqrt{b(t)-a(t)}J_0Q\right)x(t,h(t,\hat{s}))+\nonumber\\
    &\quad \sqrt{b(t)-a(t)}\left(\partial_th(t,\hat{s}) + J_1Q\right)(\partial_\sigma x(t,\sigma))\mid_{\sigma=h(t,\hat{s})},\nonumber\\
    &\stackrel{\eqref{eq:dsxhat},\eqref{eq:xhat}}{=} \left(\frac{1}{2}\frac{\dot{b}(t)-\dot{a}(t)}{b(t)-a(t)} + \frac{\dot{a}(t)+(\dot{b}(t)-\dot{a}(t))\hat{s}}{b(t)-a(t)}\partial_{\hat{s}}\right)\hat{x}(t,\hat{s})\nonumber\\
    &\quad + \hat{\mathcal{J}}(t)\hat{e}(t,\hat{s}),\nonumber\\
    &\quad= -\frac{1}{2}\frac{\dot{b}(t)-\dot{a}(t)}{b(t)-a(t)}\hat{x}(t,\hat{s}) + \partial_{\hat{s}}\left(\frac{\dot{a}+(\dot{b}-\dot{a})\hat{s}}{b-a}\hat{x}(t,\hat{s})\right)\nonumber\\
    &\quad+ \hat{\mathcal{J}}(t)\hat{e}(t,\hat{s}),\label{eq:dyn2}
\end{align}
where $\hat{\mathcal{J}}(t)\coloneqq J_0+(1/(b(t)-a(t)))J_1\partial_{\hat{s}}$. We observe that~\eqref{eq:dyn2} is exactly~\eqref{eq:dyn}, which completes this proof.

\section{Proof of Lemma~\ref{lem:power-balance}}\label{app:power-balance}
Consider the linear boundary pH system~\eqref{eq:pH-DPS} with spatial domain $\mathbb{S}_{ab}(t)$ and the Hamiltonian $\hat{H}$ in~\eqref{eq:Hhat}. By Lemma~\ref{lem:stokes-dirac} and using the fact that $J_0=-J_0^\top$, it holds, for all $\gamma(t,s)=\left(\partial_t x(t,s),f_{\partial}(t),Qx(t,s),e_{\partial}(t)\right)\in\mathcal{D}_{\mathcal{J}}(t)$, $t\in\mathbb{T}$, that
\begin{align*}
    0 &=\left\llangle\gamma(t,s),\gamma(t,s)\right\rrangle_{\mathcal{B}},\\
    &= \int_{a(t)}^{b(t)} e^\top(t,s)\partial_tx(t,s)\,\mathrm{d}s - e_\partial^\top(t) f_\partial(t),\\
    &= \int_{a(t)}^{b(t)} e^\top(t,s)\mathcal{J}e(t,s)\,\mathrm{d}s- e_\partial^\top(t) f_\partial(t),\\
    &= \int_{a(t)}^{b(t)} e^\top(t,s)J_1\partial_se(t,s)\,\mathrm{d}s - e_\partial^\top(t) f_\partial(t),
\end{align*}
where $S_1=S_1^\top\in\mathbb{R}^{r\times r}$ and $M\in\mathbb{R}^{r\times n}$ are the matrix introduced in Lemma~\ref{lem:stokes-dirac}, and
\begin{equation*}
    \begin{bmatrix}
        f_{\partial}(t)\\
        e_{\partial}(t)
    \end{bmatrix} = \begin{bmatrix}
        S_1 & -S_1\\
        I_r & I_r
    \end{bmatrix}\tr\left(Me(t,s)\right).
\end{equation*}
Applying the change of spatial variable~\eqref{eq:coc} and using~\eqref{eq:dsxhat}, we obtain
\begin{align}\label{eq:static-identity}
    &\int_0^1 \hat{e}^\top(t,\hat{s})\mathcal{\hat{J}}(t)\hat{e}(t,\hat{s})\,\mathrm{d}\hat{s} =\\
    &\quad-\left[\frac{\hat{e}^\top(t,\hat{s})M^\top S_1M\hat{e}(t,\hat{s})}{b(t)-a(t)}\right]_0^1.\nonumber
\end{align}

Differentiating $\hat{H}(\hat{x},t)$ with respect to time yields
\begin{align}
    &\frac{\mathrm{d}\hat{H}(\hat{x}(t,\cdot),t)}{\mathrm{d}t} = \int_0^1 \hat{e}^\top(t,\hat{s})\partial_t\hat{x}(t,\hat{s})\,\mathrm{d}\hat{s} \label{eq:distributed-power}\\
    &\stackrel{\eqref{eq:dyn}}{=} \int_0^1 \hat{e}^\top(t,\hat{s})\hat{\mathcal{J}}(t)\hat{e}(t,\hat{s}) -\frac{1}{2}\frac{\dot{b}(t)-\dot{a}(t)}{b(t)-a(t)}\hat{e}^\top(t,\hat{s})\hat{x}(t,\hat{s})+\nonumber\\
    &\quad \hat{e}^\top(t,\hat{s})\partial_{\hat{s}}\left(\frac{\dot{a}(t)+(\dot{b}(t)-\dot{a}(t))\hat{s}}{b(t)-a(t)}\hat{x}(t,\hat{s})\right)\,\mathrm{d}\hat{s}.\nonumber\\
    &\stackrel{\eqref{eq:static-identity}}{=} -\left[\frac{\hat{e}^\top(t,\hat{s})M^\top S_1M\hat{e}(t,\hat{s})}{b(t)-a(t)}\right]_0^1+\nonumber\\
    &\quad\frac{1}{2}\int_0^1 \partial_{\hat{s}}\left(\frac{\dot{a}(t)+(\dot{b}(t)-\dot{a}(t))\hat{s}}{b(t)-a(t)}\hat{e}^\top(t,\hat{s})\hat{x}(t,\hat{s})\right)\,\mathrm{d}\hat{s},\nonumber\\
    &= \left[\frac{1}{2}\frac{\dot{a}(t)+(\dot{b}(t)-\dot{a}(t))\hat{s}}{b(t)-a(t)}\hat{e}^\top(t,\hat{s})\hat{x}(t,\hat{s})\right]_0^1 - \nonumber\\
    &\quad\left[\frac{\hat{e}^\top(t,\hat{s})M^\top S_1M\hat{e}(t,\hat{s})}{b(t)-a(t)}\right]_0^1.\label{eq:boundary-power}
\end{align}
The proof is completed by combining~\eqref{eq:distributed-power} and~\eqref{eq:boundary-power}

\section{Proof of Theorem~\ref{thm:novel-boundary-stokes-dirac}}\label{app:novel-boundary-stokes-dirac}
Based on Definition~\ref{def:Dirac}, $\hat{\mathcal{D}}(t)$ is a Dirac structure for all $t\in\mathbb{T}$, if $\hat{\mathcal{D}}(t)=\hat{\mathcal{D}}^{\perp}(t)$ for all $t\in\mathbb{T}$, i.e., if $\hat{\mathcal{D}}(t)\subset\hat{\mathcal{D}}^{\perp}(t)$ and $\hat{\mathcal{D}}^{\perp}(t)\subset\hat{\mathcal{D}}(t)$ for all $t\in\mathbb{T}$. Before showing this, however, we observe that, by Lemma~\ref{lem:stokes-dirac} and using the fact that $J_0=-J_0^\top$, it holds, for any $(f,f_\partial,e,e_{\partial}),(\phi,\phi_{\partial},\eta,\eta_{\partial})\in\mathcal{D}_{\mathcal{J}}(t)$, $t\in\mathbb{T}$, that
\begin{align*}
    0 &= \llangle (f,f_\partial,e,e_{\partial}),(\phi,\phi_{\partial},\eta,\eta_{\partial})\rrangle_{\mathcal{B}},\\
    &= \int_{a(t)}^{b(t)} \eta^\top(t,s) \mathcal{J}e(t,s) + e^\top(t,s) \mathcal{J}\eta(t,s)\,\mathrm{d}s -\\
    &\quad f_{\partial}^\top(t)\eta_{\partial}(t) - \phi_{\partial}^\top(t) e_{\partial}(t),\\
    &= \int_{a(t)}^{b(t)} \eta^\top(t,s) J_1\partial_se(t,s) + e^\top(t,s) J_1\partial_s\eta(t,s)\,\mathrm{d}s - \\
    &\quad \phi_{\partial}^\top(t) e_{\partial}(t),\\
    &= \left[\eta^\top(t,s) J_1e(t,s)\right]_{a(t)}^{b(t)} - \phi_{\partial}^\top(t) e_{\partial}(t).
\end{align*}
It follows that, for all $e,\eta\in \hat{\mathcal{E}}(t)$, $t\in\mathbb{T}$,
\begin{align}
    &\left[\eta^\top(t,s) J_1e(t,s)\right]_{a(t)}^{b(t)} = \left(e(t,a(t))-e(t,b(t))\right)^\top\times\nonumber\\
    & M^\top S_1M\left(\eta(t,a(t))+\eta(t,b(t))\right) +\nonumber\\
    &\left(\eta(t,a(t))-\eta(t,b(t))\right)^\top \times\nonumber\\
    &M^\top S_1M\left(e(t,a(t))+e(t,b(t))\right).\label{eq:etaJe}
\end{align}

\subsection{$\hat{\mathcal{D}}(t)\subset\hat{\mathcal{D}}^\perp(t)$}
Let $\gamma=(\hat{f},\hat{f}_{\partial},\hat{e},\hat{e}_{\partial})\in\hat{\mathcal{D}}(t)$ and $\gamma_1=(\hat{\phi},\hat{\phi}_{\partial},\hat{\eta},\hat{\eta}_{\partial})\in\hat{\mathcal{D}}(t)$, $t\in\mathbb{T}$. Then, we obtain 
\begin{align}
    &\llangle \gamma(t,\hat{s}),\gamma_1(t,\hat{s})\rrangle_{\hat{\mathcal{B}}} =\int_0^1 \hat{\phi}^\top(t,\hat{s}) \hat{e}(t,\hat{s}) + \nonumber\\
    &\hat{f}^\top(t,\hat{s})\hat{\eta}(t,\hat{s}) \,\mathrm{d}\hat{s}-\hat{\phi}_{\partial}^\hop(t)\hat{e}_{\partial}(t) - \hat{f}_{\partial}^\hop(t)\hat{\eta}_{\partial}(t)=\nonumber\\
    &\int_0^1 \hat{\eta}^\top(t,\hat{s})\hat{\mathcal{J}}(t)\hat{e}(t,\hat{s}) + \hat{e}^\top(t,\hat{s})\hat{\mathcal{J}}(t)\hat{\eta}(t,\hat{s})-\nonumber\\
    &\frac{\dot{b}(t)-\dot{a}(t)}{b(t)-a(t)}\hat{\eta}^\top(t,\hat{s})Q^{-1}\hat{e}(t,\hat{s}) + \nonumber\\
    &\hat{\eta}^\top(t,\hat{s})\partial_{\hat{s}}\left(\frac{\dot{a}(t)+(\dot{b}(t)-\dot{a}(t))\hat{s}}{b(t)-a(t)}Q^{-1}\hat{e}(t,\hat{s})\right) + \nonumber\\
    &\hat{e}^\top(t,\hat{s})\partial_{\hat{s}}\left(\frac{\dot{a}(t)+(\dot{b}(t)-\dot{a}(t))\hat{s}}{b(t)-a(t)}Q^{-1}\hat{\eta}(t,\hat{s})\right)\,\mathrm{d}\hat{s} -\nonumber\\
    &\hat{\phi}_{\partial}^\hop(t)\hat{e}_{\partial}(t) - \hat{f}_{\partial}^\hop(t)\hat{\eta}_{\partial}(t).\label{eq:all-terms}
\end{align}
We proceed by considering the terms due to the internal dynamics and the time-varying boundary separately. Firstly, we find for the terms related to the moving boundaries that
\begin{align}    
    &\int_0^1 \hat{e}^\top(t,\hat{s})\partial_{\hat{s}}\left(\frac{\dot{a}(t)+(\dot{b}(t)-\dot{a}(t))\hat{s}}{b(t)-a(t)}Q^{-1}\hat{\eta}(t,\hat{s})\right) + \nonumber\\
    &\hat{\eta}^\top(t,\hat{s}) \partial_{\hat{s}}\left(\frac{\dot{a}(t)+(\dot{b}(t)-\dot{a}(t))\hat{s}}{b(t)-a(t)}Q^{-1}\hat{e}(t,\hat{s})\right)-\nonumber\\
    &\frac{\dot{b}(t)-\dot{a}(t)}{b(t)-a(t)}\eta^\top(t,\hat{s})Q^{-1}\hat{e}(t,\hat{s})\,\mathrm{d}\hat{s}=\nonumber\\
    &\int_0^1 \partial_{\hat{s}}\left(\frac{\dot{a}(t)+(\dot{b}(t)-\dot{a}(t)\hat{s}}{b(t)-a(t)}\hat{e}^\top(t,\hat{s}) Q^{-1}\hat{\eta}(t,\hat{s})\right)\,\mathrm{d}\hat{s}=\nonumber\\
    &\left[\frac{\dot{a}(t)+(\dot{b}(t)-\dot{a}(t))\hat{s}}{b(t)-a(t)}\hat{e}^\top(t,\hat{s}) Q^{-1}\hat{\eta}(t,\hat{s})\right]_0^1=\nonumber\\
    &\hat{\phi}_{\partial,2}^\hop(t)\hat{e}_{\partial,2}(t) + \hat{f}_{\partial,2}^\hop(t)\hat{\eta}_{\partial,2}(t),\label{eq:dynamic-boundary-terms}
\end{align}
where the last equality follows from the fact that, due to Assumption~\ref{asm:bounds}.\ref{item:same-sign}, $\sqrt{\dot{a}(t)}\sqrt{\dot{b}(t)}=\overline{\sqrt{\dot{a}(t)}}\overline{\sqrt{\dot{b}(t)}}$ for all $t\in\mathbb{T}$, and where we have partitioned $\hat{\phi}_{\partial}(t) = (\hat{\phi}_{\partial,1}(t),\hat{\phi}_{\partial,2}(t))$, $\hat{\eta}_{\partial}(t)=(\hat{\eta}_{\partial,1}(t),\hat{\eta}_{\partial,2}(t))$, $\hat{f}_{\partial}(t)=(\hat{f}_{\partial,1}(t),\hat{f}_{\partial,2}(t))$, and $\hat{e}_{\partial}(t)=(\hat{e}_{\partial,1}(t),\hat{e}_{\partial,2}(t))$. 
Secondly, we consider the internal dynamics, for which we find, using the fact that $J_0=-J_0^\top$, that
\begin{align}
    &\int_0^1 \hat{\eta}^\top(t,\hat{s})\hat{\mathcal{J}}(t)\hat{e}(t,\hat{s}) + \hat{e}^\top(t,\hat{s})\hat{\mathcal{J}}(t)\hat{\eta}(t,\hat{s})\,\mathrm{d}\hat{s} = \nonumber\\
    &\int_0^1 \frac{\hat{\eta}^\top(t,\hat{s}) J_1\partial_{\hat{s}}\hat{e}(t,\hat{s}) + \hat{e}^\top(t,\hat{s}) J_1\partial_{\hat{s}}\hat{\eta}(t,\hat{s})}{b(t)-a(t)}\,\mathrm{d}\hat{s},\nonumber\\
    &= \int_0^1 \frac{\partial_{\hat{s}}(\hat{\eta}^\top(t,\hat{s}) J_1\hat{e}(t,\hat{s}))}{b(t)-a(t)}\,\mathrm{d}\hat{s} = \left[\frac{\hat{\eta}^\top(t,\hat{s}) J_1\hat{e}(t,\hat{s})}{b(t)-a(t)}\right]_0^1,\nonumber\\
    &\stackrel{\eqref{eq:etaJe}}{=} \frac{\left(\hat{e}(t,0)-\hat{e}(t,1)\right)^\top M^\top S_1M\left(\hat{\eta}(t,0)+\hat{\eta}(t,1)\right)}{b(t)-a(t)} +\nonumber\\
    &\frac{\left(\hat{\eta}(t,0)-\hat{\eta}(t,1)\right)^\top M^\top S_1M \left(\hat{e}(t,0)+\hat{e}(t,1)\right)}{b(t)-a(t)}=\nonumber\\
    &\hat{\phi}^\hop_{\partial,1}(t)\hat{e}_{\partial,1}(t) + \hat{f}_{\partial,1}^{\hop}(t)\hat{\eta}_{\partial,1}(t).\label{eq:internal-dyn-terms}
\end{align}
Substituting~\eqref{eq:dynamic-boundary-terms} and~\eqref{eq:internal-dyn-terms} into~\eqref{eq:all-terms} yields
\begin{align*}
    &\llangle \gamma(t,\hat{s}),\gamma_1(t,\hat{s})\rrangle_{\hat{\mathcal{B}}} = \hat{\phi}_{\partial}^\hop(t)\hat{e}_{\partial}(t) + \hat{f}_{\partial}^\hop(t)\hat{\eta}_{\partial}(t) -\\
    &\quad \hat{\phi}_{\partial}^\hop(t)\hat{e}_{\partial}(t) - \hat{f}_{\partial}^\hop(t)\hat{\eta}_{\partial}(t) = 0.
\end{align*}
It follows that $\gamma_1\in \hat{\mathcal{D}}^\perp(t)$ and, thus, $\hat{\mathcal{D}}(t)\subset\hat{\mathcal{D}}^\perp(t)$ for all $t\in\mathbb{T}$.

\subsection{$\hat{\mathcal{D}}^\perp(t)\subset\hat{\mathcal{D}}(t)$}
Let $\gamma=(\hat{f},\hat{f}_{\partial},\hat{e},\hat{e}_{\partial})\in\hat{\mathcal{D}}(t)$ and $\gamma_1=(\hat{\phi},\hat{\phi}_{\partial},\hat{\eta},\hat{\eta}_{\partial})\in\hat{\mathcal{D}}^\perp(t)$, $t\in\mathbb{T}$. Then, 
\begin{align}
    & 0 = \llangle \gamma(t,\hat{s}),\gamma_1(t,\hat{s})\rrangle_{\hat{\mathcal{B}}}=\int_0^1 \hat{e}^\top(t,\hat{s})\hat{\phi}(t,\hat{s}) + \nonumber\\
    &\quad \hat{\eta}^\top(t,\hat{s})\hat{f}(t,\hat{s})\,\mathrm{d}\hat{s} - \hat{e}_{\partial}^\hop(t)\hat{\phi}_{\partial}(t) - \hat{f}_{\partial}^\hop(t)\hat{\eta}_{\partial}(t),\nonumber\\
    &= \int_0^1 -\frac{1}{2}\frac{\dot{b}(t)-\dot{a}(t)}{b(t)-a(t)}\hat{\eta}^\top(t,\hat{s})Q^{-1}\hat{e}(t,\hat{s})+\nonumber\\
    &\quad \hat{\eta}^\top(t,\hat{s})\partial_{\hat{s}}\left(\frac{\dot{a}(t)+(\dot{b}(t)-\dot{a}(t))\hat{s}}{b(t)-a(t)}Q^{-1}\hat{e}(t,\hat{s})\right)+\nonumber\\
    &\quad \hat{\eta}^\top(t,\hat{s})\hat{\mathcal{J}}(t)\hat{e}(t,\hat{s})+\hat{\phi}^\top(t,\hat{s})\hat{e}(t,\hat{s})\,\mathrm{d}\hat{s} -\nonumber\\
    &\quad \hat{e}^\hop_{\partial}(t)\hat{\phi}_{\partial}(t)-\hat{f}_{\partial}^\hop(t)\hat{\eta}_{\partial}(t).\label{eq:all-terms-2}
\end{align}
As before, let us treat the terms due to the internal dynamics and those due to the moving boundaries separately. First, we consider terms related to the internal dynamics, which we manipulate as follows:
\begin{align}
    &\int_0^1 \hat{\eta}^\top(t,\hat{s})\hat{\mathcal{J}}(t)\hat{e}(t,\hat{s})\,\mathrm{d}\hat{s} = \nonumber\\
    &\int_0^1 \hat{\eta}^\top(t,\hat{s})\left(J_0 + \frac{1}{b(t)-a(t)}J_1\partial_{\hat{s}}\right)\hat{e}(t,\hat{s})\,\mathrm{d}\hat{s},\nonumber\\
    &= \int_0^1 \frac{\partial_{\hat{s}}(\hat{e}^\top(t,\hat{s}) J_1\hat{\eta}(t,\hat{s}))}{b(t)-a(t)} - \\
    &\quad \hat{e}^\top(t,\hat{s})\left(J_0 + \frac{1}{b(t)-a(t)}J_1\partial_{\hat{s}}\right)\hat{\eta}(t,\hat{s})\,\mathrm{d}\hat{s},\nonumber\\
    &= \left[\frac{\hat{e}^\top(t,\hat{s}) J_1\hat{\eta}(t,\hat{s})}{b(t)-a(t)}\right]_0^1 - \int_0^1 \hat{e}^\top(t,\hat{s}) \hat{\mathcal{J}}(t)\hat{\eta}(t,\hat{s})\,\mathrm{d}\hat{s},\nonumber\\
    &\stackrel{\eqref{eq:etaJe}}{=} \frac{\left(\hat{e}(t,0)-\hat{e}(t,1)\right)^\top M^\top S_1M\left(\hat{\eta}(t,0)+\hat{\eta}(t,1)\right)}{b(t)-a(t)} +\nonumber\\
    &\frac{\left(\hat{\eta}(t,0)-\hat{\eta}(t,1)\right)^\top M^\top S_1M \left(\hat{e}(t,0)+\hat{e}(t,1)\right)}{b(t)-a(t)}-\nonumber\\
    &\int_0^1 \hat{e}^\top(t,\hat{s}) \hat{\mathcal{J}}(t)\hat{\eta}(t,\hat{s})\,\mathrm{d}\hat{s}.\label{eq:internal-dyn-terms-2}
\end{align}
Secondly, we consider the terms due to the moving boundaries, which we manipulate as follows:
\begin{align}
    &\int_0^1 -\frac{1}{2}\frac{\dot{b}(t)-\dot{a}(t)}{b(t)-a(t)}\hat{\eta}^\top(t,\hat{s})Q^{-1}\hat{e}(t,\hat{s}) +\nonumber\\
    &\quad \hat{\eta}^\top(t,\hat{s})\partial_{\hat{s}}\left(\frac{\dot{a}(t)+(\dot{b}(t)-\dot{a}(t))\hat{s}}{b(t)-a(t)}Q^{-1}\hat{e}(t,\hat{s})\right)\,\mathrm{d}\hat{s}\nonumber\\
    &=\int_0^1 \partial_{\hat{s}}\left(\frac{\dot{a}(t)+(\dot{b}(t)-\dot{a}(t))\hat{s}}{b(t)-a(t)}\hat{\eta}^\top(t,\hat{s}) Q^{-1}\hat{e}(t,\hat{s})\right) +\nonumber\\
    &\quad \frac{1}{2}\frac{\dot{b}(t)-\dot{a}(t)}{b(t)-a(t)}\hat{e}^\top(t,\hat{s}) Q^{-1}\hat{\eta}(t,\hat{s}) -\nonumber\\
    &\quad \hat{e}^\top(t,\hat{s})\partial_{\hat{s}}\left(\frac{\dot{a}(t)+(\dot{b}(t)-\dot{a}(t))\hat{s}}{b(t)-a(t)} Q^{-1}\hat{\eta}(t,\hat{s})\right)\,\mathrm{d}\hat{s} = \nonumber\\
    &-\int_0^1 \hat{e}^\top(t,\hat{s}) Q^{-1}\partial_{\hat{s}}\left(\frac{\dot{a}(t)+(\dot{b}(t)-\dot{a}(t))\hat{s}}{b(t)-a(t)}\hat{\eta}(t,\hat{s})\right)\nonumber\\
    &\quad -\frac{1}{2}\frac{\dot{b}(t)-\dot{a}(t)}{b(t)-a(t)}\hat{e}^\top(t,\hat{s})Q^{-1}\hat{\eta}(t,\hat{s})\,\mathrm{d}\hat{s}\nonumber\\
    &+\left[\frac{\dot{a}(t)+(\dot{b}(t)-\dot{a}(t))\hat{s}}{b(t)-a(t)}\hat{e}^\top(t,\hat{s})Q^{-1}\hat{\eta}(t,\hat{s})\right]_0^1.\label{eq:dynamic-boundary-terms-2}
\end{align}
Substituting~\eqref{eq:internal-dyn-terms-2} and~\eqref{eq:dynamic-boundary-terms-2} into~\eqref{eq:all-terms-2} together with the expressions for $\hat{e}_{\partial}(t)$ and $\hat{f}_{\partial}(t)$ in~\eqref{eq:boundary-ports} allows us to factor as follows:
\begin{align}
    0 &= \llangle \gamma(t,\hat{s}),\gamma_1(t,\hat{s})\rrangle_{\hat{\mathcal{B}}}=\int_0^1 \hat{e}^\top(t,\hat{s})\left(\hat{\phi}(t,\hat{s})-  \right.\label{eq:all-terms-final}\\
    &\left.\hat{\mathcal{J}}(t)\hat{\eta}(t,\hat{s}) +\frac{1}{2}\frac{\dot{b}(t)-\dot{a}(t)}{b(t)-a(t)}Q^{-1}\hat{\eta}(t,\hat{s})-\right.\nonumber\\
    &\quad\left.\partial_{\hat{s}}\left(\frac{\dot{a}(t)+(\dot{b}(t)-\dot{a}(t))\hat{s}}{b(t)-a(t)}Q^{-1}\hat{\eta}(t,\hat{s})\right)\right)\mathrm{d}\hat{s}+\nonumber\\
    &\frac{\left(\hat{e}(t,0)-\hat{e}(t,1)\right)^\top M^\top S_1}{b(t)-a(t)}\times\nonumber\\
    &\quad\left(M\left(\hat{\eta}(t,0)+\hat{\eta}(t,1)\right)-\sqrt{b(t)-a(t)}\hat{\eta}_{\partial,1}(t)\right)+\nonumber\\
    &\frac{\left(\hat{e}(t,0)+\hat{e}(t,1)\right)^\top M^\top}{b(t)-a(t)}\times\nonumber\\
    &~\left(S_1M\left(\hat{\eta}(t,0)-\hat{\eta}(t,1)\right)-\sqrt{b(t)-a(t)}\hat{\phi}_{\partial,1}(t)\right)-\nonumber\\
    &\left(\frac{\sqrt{\dot{a}(t)}Q^{-1}\hat{e}(t,0)+\sqrt{\dot{b}(t)}Q^{-1}\hat{e}(t,1)}{2\sqrt{b(t)-a(t)}}\right)^\hop\hat{\phi}_{\partial,2}(t)\nonumber\\
    &+\left(\frac{\overline{\sqrt{\dot{a}(t)}}\hat{e}(t,0)-\overline{\sqrt{\dot{b}(t)}}\hat{e}(t,1)}{\sqrt{b(t)-a(t)}}\right)^\hop\hat{\eta}_{\partial,2}(t)+\nonumber\\
    &\left[\frac{\dot{a}(t)+(\dot{b}(t)-\dot{a}(t))\hat{s}}{b(t)-a(t)}\hat{e}^\top(t,\hat{s})Q^{-1}\hat{\eta}(t,\hat{s})\right]_0^1.\nonumber
\end{align}
where we have partitioned $\hat{\phi}_{\partial}(t) = (\hat{\phi}_{\partial,1}(t),\hat{\phi}_{\partial,2}(t))$ and $\hat{\eta}_{\partial}(t)=(\hat{\eta}_{\partial,1}(t),\hat{\eta}_{\partial,2}(t))$. Since~\eqref{eq:all-terms-final} should hold for all $\hat{e}\in \hat{\mathcal{E}}(t)$, $t\in\mathbb{T}$, and for all $a,b\in C(\mathbb{T},\mathbb{R})$ satisfying Assumption~\ref{asm:bounds} (so also for $\dot{a}=\dot{b}=0$) and since $M$ is full rank, it follows that
\begin{align*}
    \hat{\phi}(t,\hat{s}) &= -\frac{1}{2}\frac{\dot{b}(t)-\dot{a}(t)}{b(t)-a(t)}Q^{-1}\hat{\eta}(t,\hat{s}) + \hat{\mathcal{J}}(t)\hat{\eta}(t,\hat{s})+\\
    &\quad \partial_{\hat{s}}\left(\frac{\dot{a}(t)+(\dot{b}(t)-\dot{a}(t))\hat{s}}{b(t)-a(t)}Q^{-1}\hat{\eta}(t,\hat{s})\right),\\
    \hat{\phi}_{\partial,1}(t) &= \frac{1}{\sqrt{b(t)-a(t)}}S_1M(\hat{\eta}(t,0)-\hat{\eta}(t,1)),\\
    \hat{\eta}_{\partial,1}(t) &= \frac{1}{\sqrt{b(t)-a(t)}}M(\hat{\eta}(t,0)+\hat{\eta}(t,1)),
\end{align*}
which are precisely the conditions in the definition of $\hat{\mathcal{D}}(t)$ in Theorem~\ref{thm:novel-boundary-stokes-dirac}. Substitution in~\eqref{eq:all-terms-final} leaves us with 
\begin{align}
    0 &= -\left(\frac{\sqrt{\dot{a}(t)}Q^{-1}\hat{e}(t,0)+\sqrt{\dot{b}(t)}Q^{-1}\hat{e}(t,1)}{2\sqrt{b(t)-a(t)}}\right)^\hop\hat{\phi}_{\partial,2}(t)\nonumber\\
    &+\left(\frac{\overline{\sqrt{\dot{a}(t)}}\hat{e}(t,0)-\overline{\sqrt{\dot{b}(t)}}\hat{e}(t,1)}{\sqrt{b(t)-a(t)}}\right)^\hop\hat{\eta}_{\partial,2}(t)+\nonumber\\
    &\left[\frac{\dot{a}(t)+(\dot{b}(t)-\dot{a}(t))\hat{s}}{b(t)-a(t)}\hat{e}^\top(t,\hat{s})Q^{-1}\hat{\eta}(t,\hat{s})\right]_0^1.\label{eq:residual}
\end{align}
Using Assumption~\ref{asm:bounds}.\ref{item:same-sign}, whereby $\sqrt{\dot{a}(t)\dot{b}(t)}=\overline{\sqrt{\dot{a}(t)\dot{b}(t)}}$ for all $t\in\mathbb{T}$, we can factor
\begin{align}
    &\left[\frac{\dot{a}(t)+(\dot{b}(t)-\dot{a}(t))\hat{s}}{b(t)-a(t)}\hat{e}^\top(t,\hat{s})Q^{-1}\hat{\eta}(t,\hat{s})\right]_0^1 = \label{eq:factored2}\\
    &-\frac{1}{2}\left(\frac{\sqrt{\dot{a}(t)}Q^{-1}\hat{e}(t,0)+\sqrt{\dot{b}(t)}Q^{-1}\hat{e}(t,1)}{\sqrt{b(t)-a(t)}}\right)^\hop\times\nonumber\\
    &\quad\left(\frac{\overline{\sqrt{\dot{a}(t)}}\hat{\eta}(t,0)-\overline{\sqrt{\dot{b}(t)}}\hat{\eta}(t,1)}{\sqrt{b(t)-a(t)}}\right) - \nonumber\\
    &\left(\frac{\overline{\sqrt{\dot{a}(t)}}\hat{e}(t,0)-\overline{\sqrt{\dot{b}(t)}}\hat{e}(t,1)}{\sqrt{b(t)-a(t)}}\right)^\hop\times\nonumber\\
    &\quad \left(\frac{1}{2}Q^{-1}\frac{\sqrt{\dot{a}(t)}\hat{\eta}(t,0)+\sqrt{\dot{b}(t)}\hat{\eta}(t,1)}{\sqrt{b(t)-a(t)}}\right).\nonumber
\end{align}
Substituting~\eqref{eq:factored2} into~\eqref{eq:residual}, we find that
\begin{align}
    0 &= -\frac{1}{2}\left(\frac{\sqrt{\dot{a}(t)}Q^{-1}\hat{e}(t,0)+\sqrt{\dot{b}(t)}Q^{-1}\hat{e}(t,1)}{\sqrt{b(t)-a(t)}}\right)^\hop\times\nonumber\\
    &\quad\left(\hat{\phi}_{\partial,2}(t)+\frac{\overline{\sqrt{\dot{a}(t)}}\hat{\eta}(t,0)-\overline{\sqrt{\dot{b}(t)}}\hat{\eta}(t,1)}{\sqrt{b(t)-a(t)}}\right) - \nonumber\\
    &\left(\frac{\overline{\sqrt{\dot{a}(t)}}\hat{e}(t,0)-\overline{\sqrt{\dot{b}(t)}}\hat{e}(t,1)}{\sqrt{b(t)-a(t)}}\right)^\hop\times\nonumber\\
    &\quad \left(-\hat{\eta}_{\partial,2}(t)+\frac{1}{2}Q^{-1}\frac{\sqrt{\dot{a}(t)}\hat{\eta}(t,0)+\sqrt{\dot{b}(t)}\hat{\eta}(t,1)}{\sqrt{b(t)-a(t)}}\right),\label{eq:final-final}
\end{align}
for all $\hat{e}\in\hat{\mathcal{E}}(t)$, $t\in\mathbb{T}$, and $a,b\in C(\mathbb{T},\mathbb{R})$ satisfying Assumption~\ref{asm:bounds}. 
It follows from~\eqref{eq:final-final} that 
\begin{equation*}
    \begin{bmatrix}
        \hat{\phi}_{\partial,2}(t)\\
        \hat{\eta}_{\partial,2}(t)
    \end{bmatrix} = \frac{1}{\sqrt{b(t)-a(t)}}\begin{bmatrix}
        -I_n & I_n\\
        \frac{1}{2}Q^{-1} & \frac{1}{2}Q^{-1}
    \end{bmatrix}\begin{bmatrix}
        \sqrt{\dot{a}(t)}\hat{\eta}(t,0)\\
        \sqrt{\dot{b}(t)}\hat{\eta}(t,1)        
    \end{bmatrix},
\end{equation*}
 Thus, it holds that $\gamma_1\in\hat{\mathcal{D}}(t)$. Hence, we have shown that $\hat{\mathcal{D}}(t)\subset\hat{\mathcal{D}}^{\perp}(t)$ for all $t\in\mathbb{T}$, which completes the proof.
\end{document}